\documentclass[preprint,12pt]{elsarticle}

\usepackage[english]{babel}

\usepackage[utf8]{inputenc}

\usepackage{color}

\usepackage{hhline}
\usepackage{multirow}
\usepackage{tabularx}
\usepackage{hyperref}
\usepackage{makecell}
\usepackage{amsmath,amssymb}
\usepackage{amsthm}
\usepackage{siunitx}
\usepackage{placeins}
\usepackage[T1]{fontenc}
\usepackage{subfig}
\usepackage{longtable}
\usepackage{graphicx}
\usepackage{tikz}

\usepackage[left=1in, right=1in, top=1in, bottom=1in]{geometry}
\usepackage{setspace}
\newcommand{\cE}{{\cal E}}

\newcommand{\cV}{{\cal V}}

\newcommand{\pdiff}[2]{\dfrac{\partial{#1}}{\partial{#2}}}

\newcommand{\diff}[2]{\dfrac{d{#1}}{d{#2}}}
\newcommand{\pd}[1]{\partial_{#1}}
\newcommand{\pdt}{\partial_t}
\newcommand{\pdx}{\partial_x}
\newcommand{\bv}[1]{{\bf #1}}

\newcommand{\IR}{{\mathbb{R}}}

\newcommand{\av}[1]{\bar{\bv{#1}}}

\renewcommand{\leq}{\leqslant}
\renewcommand{\geq}{\geqslant}

\newcommand{\const}{\textrm{const}}

\begin{document}

\begin{frontmatter}



\title{Stochastic Finite Volume Method for Uncertainty Quantification of Transient Flow in Gas Pipeline Networks}


\author[LANLT5]{S. Tokareva}

\author[LANLT5]{A. Zlotnik}

\author[LANLT5]{V. Gyrya}

\affiliation[LANLT5]{organization={Los Alamos National Laboratory},
            city={Los Alamos},
            state={NM},
            country={USA}}

\begin{abstract}
We develop a weakly intrusive framework to simulate the propagation of uncertainty in solutions of generic hyperbolic partial differential equation systems on graph-connected domains with nodal coupling and boundary conditions. The method is based on the Stochastic Finite Volume (SFV) approach, and can be applied for uncertainty quantification (UQ) of the dynamical state of fluid flow over actuated transport networks.  The numerical scheme has specific advantages for modeling intertemporal uncertainty in time-varying boundary parameters, which cannot be characterized by strict upper and lower (interval) bounds.  We describe the scheme for a single pipe, and then formulate the controlled junction Riemann problem (JRP) that enables the extension to general network structures.  We demonstrate the method's capabilities and performance characteristics using a standard benchmark test network.

\end{abstract}


\end{frontmatter}

\section{Introduction}

Many problems in physics and engineering are modeled by hyperbolic systems of conservation or balance laws. Several prominent examples include the Shallow Water Equations of hydrology \cite{fjordholm2011well}, the Euler Equations for inviscid, compressible flow \cite{jameson1985numerical}, and the Magnetohydrodynamic (MHD) equations of plasma physics \cite{beg2013numerical}. Many efficient numerical methods have been developed to approximate the solutions of systems of conservation laws \cite{Godlewski-Raviart,LeVeque}, e.g. finite volume schemes \cite{morton2007finite} or discontinuous Galerkin methods \cite{krivodonova2004shock}. The classical assumption in designing efficient numerical methods for hyperbolic systems is that the initial data, boundary conditions, and coefficients of the model are known exactly, i.e., they are deterministic. However, in many practical applications it is not always possible to obtain exact data due to, for example, measurement, prediction, or modeling errors. 

Numerous studies describe incomplete information in the uncertain data mathematically as random fields. Such data are described in terms of statistical quantities of interest such as the mean, variance, and higher statistical moments, and in some cases the distribution law of the stochastic data is also assumed to be known.  Numerical techniques have been developed to quantify process uncertainty by computing the mean flow and its statistical moments.  They are based on approaches such as the Monte-Carlo (MC) \cite{MishraSchwab2010}, multi-level Monte Carlo (MLMC) \cite{MishraSchwabSukys2011}, generalized polynomial chaos (gPC) \cite{Poette,Lin-Su-Karniadakis-1}, the probabilistic collocation method (PCM) \cite{Lin-Su-Karniadakis-2}, Godunov schemes \cite{Troyen-2}, and stochastic Galerkin (sG) projections \cite{Troyen-1,Gottlieb-Xiu}.

Uncertainty quantification methods can be roughly classified into intrusive and non-intrusive schemes. \textit{Non-intrusive} UQ methods allow the re-use of an existing deterministic code as a black box without any modifications. An example of such an approach is the Monte Carlo method: given a large number of realizations of the random parameters, we generate the corresponding outputs from our deterministic code and then process this data to obtain the statistical mean and variances. The possibility to reuse existing deterministic codes is a clear advantage of non-intrusive approaches. However, methods such as Monte Carlo tend to require a prohibitively large number of evaluations of the deterministic solutions.  \textit{Intrusive} UQ methods, on the contrary, require modifications to the algorithms and their implementations in the deterministic computational scheme. For example, a well known and widely used sG method for hyperbolic systems transforms the original PDEs for the primary variables into a set of PDEs that are defined with respect to the polynomial expansion coefficients \cite{Troyen-1}. This reduces the simulation time, but might pose mathematical challenges, such as the loss of hyperbolicity of the transformed PDEs.

The application of uncertainty quantification to hyperbolic conservation laws on graphs has found compelling motivation from the need to model and control vehicle traffic \cite{gerster2021stability} and gas pipeline flows \cite{gerster2019polynomial}.  These recent studies relied on tailored polynomial chaos expansions and the stochastic Galerkin method for forward simulation of the stochastic PDEs subject to distributional uncertainty in initial conditions and parameters  \cite{gerster2019hyperbolic}, and primarily examined the problem of stabilization.  The representation of uncertainty in gas pipeline network flows requires the solution of stochastic PDEs on a large number of one-dimensional domains connected in a graph.  Kirchhoff-like conditions are used to specify mass conservation at graph vertices, and may moreover include the effect of gas compressors that are used to actuate flow throughout the network.  Two key challenges are generalizability and scalability of the technique to arbitrary and large graphs, and handling of inter-temporal uncertainty.  Gas pipelines often experience temporary increases in gas demand of limited duration that start at an uncertain time, so that two possible boundary condition profiles would not follow strict pointwise ordering in time.  As a result, using distributional uncertainty on the boundary parameter point-wise in time would result in an excessive uncertainty estimate.  The ability to efficiently compute an accurate uncertainty quantification for system pressures and flows in such scenarios motivates our investigation.

In this study, we present a semi-intrusive Stochastic Finite Volume (SFV) method to quantify the uncertainties that arise due to random model parameters in the underlying hyperbolic PDE system, including initial conditions, uncertain constants, and complex time-dependent distributions on boundary conditions.  The SFV method requires some modifications of the deterministic code that is used for solving a standard initial boundary value problem (IBVP) for hyperbolic conservation laws, which only involve additional integration of the numerical fluxes over the cells in the stochastic space and can therefore be considered mild.  This approach preserves the hyperbolicity of the model, and at the same time is more computationally efficient than the Monte Carlo method.  To enable the generalization to graph domains, we formulate a junction Riemann problem (JRP) method for propagating uncertainty solutions through graph nodes with flow controllers. We apply the scheme for uncertainty quantification of the dynamical state of fluid flow over transport networks with nodal coupling conditions.

The remainder of the paper is organized as follows. In Section \ref{sec:sfvmethod}, we provide a general presentation of the Stochastic Finite Volume (SFV) method, followed by a brief description of the one-dimensional model of gas flow in a pipe and the underlying model uncertainty that accounts for inter-temporal factors. Section \ref{sec:networkuq} then extends the SFV methods for UQ in PDEs on graphs by introducing the JRP for the considered model of gas flow, and summarizes standard network modeling for gas pipeline systems as the setting of application. Finally, numerical results for a single pipe and a small test network are given in Section \ref{sec:computation}.

\section{Uncertainty Quantification for Gas Flow in a Pipe} \label{sec:sfvmethod}

Our goal is to develop a simple, flexible, and extensible stochastic representation of gas flow and pressure in a pipe subject to transient boundary conditions that will be easily used as a building block for more complex simulations.  In this section we describe in general the SFV method for hyperbolic conservation laws, and then the basic problem for gas flow in a pipe.  We then explain the challenge of intertemporal uncertainty modeling, which to our knowledge has not been addressed using standard methods for even this basic problem. 

\subsection{Stochastic Finite Volume method} \label{sec:sfv_general}

Here we present the Stochastic Finite Volume Method (SFVM) for uncertainty quantification for general conservation laws \cite{Abgrall,BAR1}, which is based on the finite volume framework. The SFVM is formulated to numerically solve the system of conservation laws with sources of randomness in flux coefficients, initial and boundary data.

Consider the hyperbolic system of conservation laws with random flux coefficients
\begin{equation}
\label{eq:MCL}
\pdiff{\bv{U}}{t} + \nabla_x\cdot\bv{F}(\bv{U},\omega) = \bv{0}
\end{equation}
defined over $t > 0$ on a spatial domain $\bv{x} = (x_1, x_2, x_3) \in D_x \subset \IR^3$, for $p$ conserved values $\bv{U} = [u_1,\dots,u_p]^\top$, and fluxes $\bv{F} = [\bv{F}_1, \bv{F}_2, \bv{F}_3]$, $\bv{F}_k = [f_1^k,\dots,f_p^k]^\top$ for  $k = 1,2,3$.  A stochastic IBVP is obtained by specifying random initial data
\begin{equation}
\label{eq:MCL-IC}
\bv{U}(\bv{x},0,\omega) = \bv{U}_0(\bv{x},\omega), \ \ \bv{x} \in D_x, \ \omega \in \Omega, 
\end{equation}
and possibly random boundary conditions
\begin{equation}
\label{eq:MCL-BC}
\bv{U}(\bv{x},t,\omega) = \bv{U}_B(t,\omega), \ \bv{x} \in \partial D_x, \ \omega \in \Omega.
\end{equation}

A seminal study has developed a mathematical framework of \textit{random entropy solutions} for such scalar conservation laws of form \eqref{eq:MCL}--\eqref{eq:MCL-IC}, and proved existence and uniqueness of such solutions for scalar conservation laws in multiple dimensions \cite{MishraSchwab2010}. Furthermore, the existence of statistical quantifiers of the random entropy solution, such as the statistical mean and $k$-point spatio-temporal correlation functions, have been proven to exist under suitable assumptions on the random initial data.

We parametrize all the random inputs in the equations \eqref{eq:MCL}--\eqref{eq:MCL-IC} using the random variable $\bv{y} = \bv{Y}(\omega)$, which takes values in $D_y \subset \IR^q$, and rewrite the stochastic conservation law in the parametric form:
\begin{subequations}
\begin{align}
\label{eq:MCL1}
 \pd{t} \bv{U} + \nabla_x\cdot\bv{F}(\bv{U},\bv{y}) = \bv{0}, \quad \bv{x} \in D_x \ \subset \IR^3, \ \bv{y} \in D_y\subset \IR^q, \ t > 0; \\
\label{eq:MCL-IC1}
 \bv{U}(\bv{x},0,\bv{y}) = \bv{U}_0(\bv{x},\bv{y}), \quad \bv{x} \in D_x \ \subset \IR^3, \ \bv{y} \in D_y\subset \IR^q. \qquad
\end{align}
\end{subequations}
Let us now define the two spaces
\begin{equation}
    \mathcal{T}_x = \bigcup_{i=1}^{N_x} K_x^i \quad \text{and} \quad \mathcal{C}_y = \bigcup_{j=1}^{N_y}K_y^j,
\end{equation}
such that $\mathcal{T}_x$ is the triangulation of the computational domain $D_x$ in the physical space, and $\mathcal{C}_y$ is the Cartesian grid in the domain $D_y$ of the parameterized probability space.  We further assume the existence of a probability density function $\mu(\bv{y})$ on the space $D_y$, and compute the expectation of the $n$-th solution component of the conservation law \eqref{eq:MCL1}--\eqref{eq:MCL-IC1} as
\begin{equation} \label{eq:exp_nth}
\mathbb{E}[u_n] = \int\limits_{D_y} u_n \mu(\bv{y})\, d\bv{y}, \ n = 1,\dots,p.
\end{equation}
The SFVM scheme \cite{SFVSpringer2014} can be obtained from the integral form of equations \eqref{eq:MCL1}--\eqref{eq:MCL-IC1}:
\begin{equation} \label{eq:sfv_int0}
\iint\limits_{K_y^j\,K_x^i} \pd{t} \bv{U} \,\mu(\bv{y})\, d\bv{x}d\bv{y} + \iint\limits_{K_y^j\,K_x^i} \nabla_x\cdot\bv{F}(\bv{U},\bv{y}) \,\mu(\bv{y})\, d\bv{x}d\bv{y} = \bv{0}.
\end{equation}
Introducing the cell average
\begin{equation} \label{eq:sfv_cellav}
\av{U}_{ij}(t) = \dfrac{1}{|K_x^i||K_y^j|}\iint\limits_{K_y^j\,K_x^i} \bv{U}(\bv{x},t,\bv{y})\mu(\bv{y})\,d\bv{x}d\bv{y}
\end{equation}
with the cell volumes
\begin{equation} \label{eq:sfv_cellvol}
|K_x^i| = \int\limits_{K_x^i} d\bv{x}, \quad |K_y^j| = \int\limits_{K_y^j} \mu(\bv{y})\,d\bv{y},
\end{equation}
and performing the partial integration over $K_x^i$, we obtain
\begin{equation} \label{eq:sfv_partint}
\diff{\bar{\bv{U}}_{ij}}{t} + \dfrac{1}{|K_x^i||K_y^j|}\int\limits_{K_y^j} \bigg[\int\limits_{\partial K_x^i} \bv{F}(\bv{U},\bv{y})\cdot \bv{n}\,dS \bigg] \mu(\bv{y})\,d\bv{y} = \bv{0}.
\end{equation}
Next, we use any standard numerical flux approximation, which we denote henceforth by $\hat{\bv{F}}\big(\tilde{\bv{U}}_L(\bv{x},t,\bv{y}),\tilde{\bv{U}}_R(\bv{x},t,\bv{y}),\bv{y}\big)$, to replace the discontinuous flux, $\bv{F}(\bv{U},\bv{y})\cdot \bv{n}$, through the element interface. Here the term $\tilde{\bv{U}}_{L}$ and $\tilde{\bv{U}}_{R}$ denote the boundary extrapolated solution values at the edge of the cell $K_x^i$, obtained by the high order reconstruction from the cell averages, see Section \ref{sec:tvd_recon} for details. The complete numerical flux is  approximated by a suitable quadrature rule as
\begin{equation}
\label{NumFlux}
\av{F}_{ij}(t) = \dfrac{1}{|K_y^j|} \int\limits_{K_y^j}  \bigg[\int\limits_{\partial K_x^i} \hat{\bv{F}}(\tilde{\bv{U}}_L,\tilde{\bv{U}}_R,\bv{y}) \bigg] \mu(\bv{y})\,d\bv{y} \approx \dfrac{1}{|K_y^j|} \sum\limits_{\bv{m}} \hat{\mathcal{F}}(t,\bv{y}_{\bv{m}}) \mu(\bv{y}_{\bv{m}}) w_{\bv{m}},  
\end{equation}
where $\hat{\mathcal{F}}$ denotes the flux integral over the physical cell, $\bv{m} = (m_1,\dots,m_q)$ is the multi-index, and $\bv{y}_{\bv{m}}$ and $w_{\bv{m}}$ are quadrature nodes and weights, respectively.  The SFV method then results in the solution of the following ODE system:
\begin{equation}
\label{SFV}
\diff{\bar{\bv{U}}_{ij}}{t} + \dfrac{1}{|K_x^i|} \av{F}_{ij}(t) = \bv{0}, \quad \forall i = 1,\dots,N_x, \,\, \forall j = 1,\dots,N_y.
\end{equation}
Therefore, to obtain the high-order scheme we first need to provide the high-order flux approximation based, for example, on essentially non-oscillatory (ENO) or weighted essentially non-oscillatory (WENO) reconstruction in the physical space \cite{shu99,zhang2016eno}. Second, we have to guarantee the high-order integration in \eqref{NumFlux} also by applying the ENO/WENO reconstruction in the stochastic space and choosing the suitable quadrature rule. Finally, we need the high-order time-stepping algorithm to solve the ODE system \eqref{SFV}, such as the Runge-Kutta method.

\subsection{SFV Method Applied to the 1-Dimensional Model of Gas Flow} \label{sec:sfv_1d}

Gas flow in a long pipe in the regime without waves or shocks can be sufficiently described \cite{misra2020monotonicity,osiadacz1984simulation} by the  one-dimensional PDE system 
\begin{subequations}\label{Eq:gas0}
\begin{align}
& \pdt \rho + \pdx q = 0, \\
& \pdt q + a^2\pdx \rho = -\dfrac{\lambda}{2D}\dfrac{q|q|}{\rho}.
\end{align}
\end{subequations}
Here $\rho$ and $q$ are density and (per area) mass flux, $D$ is pipe diameter, $\lambda$ is the Darcy-Weisbach friction factor, and $a$ is the wave speed. The PDEs are solved for a pipe of length $L$, so that $x \in [0,L]$, and over the time interval $t \in [0,T]$.  Specifying initial conditions $\rho(0,x)=\rho_0(x)$ and $q(0,x) = q_0(x)$, as well as boundary conditions $\rho(0,t) = s(t)$ and $q(L,t) = d(t)$, results in a well-posed IBVP.

\begin{figure}[h!]
    \centering
    \includegraphics[width=0.9\textwidth]{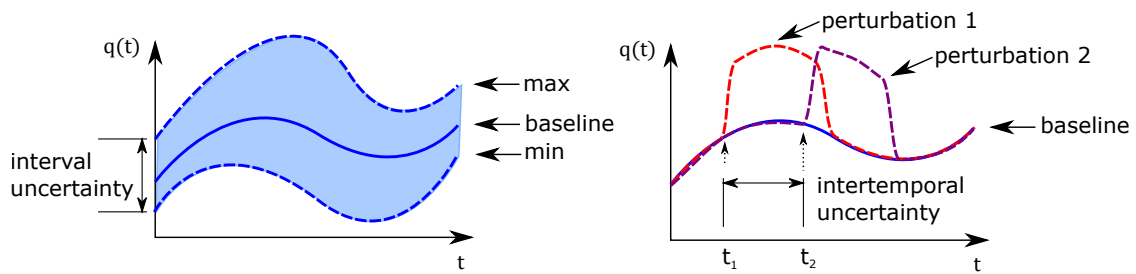}
    \caption{Left: interval uncertainty between ordered max. and min. profiles.  Right: intertemporal uncertainty of perturbations.}
    \label{fig:uncertainties-peak}
\end{figure}

We aim to construct a method that can be used to model the propagation of both \textit{interval} and \textit{intertemporal} uncertainties, as illustrated in Fig.~\ref{fig:uncertainties-peak}.  Interval uncertainty is used to denote the possibility of time-dependent parameters delimited by upper and lower bounding functions.  Interpemporal uncertainty denotes an augmentation in a time-varying parameter  starting at a time that is randomly distributed (e.g. uniformly on an interval).  An example of intertemporal uncertainty in the boundary flow $d(t)$ is
\begin{equation}
\hat{d}(t;\omega) = 
\begin{cases}
& d(t) + \delta d, \quad \text{if} \quad t \in (t_p(\omega), t_p(\omega) + \delta t_p), \\
& d(t), \quad \text{otherwise},
\end{cases}
\end{equation}
where $t_p(\omega)$ denotes a random augmentation time.  From Fig.~\ref{fig:uncertainties-peak}, we see that interval uncertainty modeling includes a much greater collection of possible time-varying parameters.

In the following, we explain in more detail the application of SFV method for a one-dimensional hyperbolic PDE system such as \eqref{Eq:gas0}. In order to simplify the presentation, we will assume that the random inputs can be represented by a single uniformly distributed random variable.  Consider the one-dimensional model of gas flow with uncertainty appearing in the initial and boundary data and source terms: 
\begin{align}
& \pdiff{\bv{U}}{t}+\pdiff{\bv{F}}{x}=\bv{S}(\bv{U},\omega), \quad x \in (x_L,x_R), \ t>0, \ \omega \in \Omega, \label{eq:1D-PDE} \\
& \bv{U}(x,0,\omega) = \bv{U}_0(x,\omega),  \label{eq:1D-IC} \\
& \bv{U}(0,t,\omega) = \bv{U}_B^0(t,\omega), \label{eq:1D-BC-0} \\
& \bv{U}(L,t,\omega) = \bv{U}_B^L(t,\omega). \label{eq:1D-BC-L}
\end{align}
For the model of gas flow in a pipe of interest here, we have 
\begin{equation}
\bv{U} = \left[\begin{array}{c} \rho \\ q \end{array}\right], \qquad \bv{F} = \left[\begin{array}{c} q \\ a^2\rho \end{array}\right],  \qquad \bv{S} =  \left[\begin{array}{c} 0 \\ -\frac{\lambda}{2D}\frac{q|q|}{\rho} \end{array}\right].
\end{equation}

Assume that the random data $\omega$ can be parametrized using one stochastic variable $y = Y (\omega) \in \IR$ with probability density function $\mu(y)$, and therefore the initial and boundary data and sources take the form
\begin{subequations}
\begin{alignat}{2}
\bv{U}(x, 0, Y(\omega)) &= \bv{U}_0(x, y), \quad &  y \in (y_L, y_R), \label{eq:1D-IC-par}\\
\bv{U}(0,t,Y(\omega)) &= \bv{U}_B^0(t,y), &  y \in (y_L, y_R), \\
\bv{U}(L,t,Y(\omega)) &= \bv{U}_B^L(t,y), &  y \in (y_L, y_R), \\
\bv{S}(\bv{U},Y(\omega)) &= \bv{S}(\bv{U},y), &  y \in (y_L, y_R). \label{eq:1D-SRC-par}
\end{alignat}
\end{subequations}
After this transformation, the parametric form of the stochastic PDE is
\begin{align}
& \pdiff{\bv{U}}{t}+\pdiff{\bv{F}}{x}=\bv{S}(\bv{U},y), \quad x \in (x_L,x_R), \ y \in (y_L, y_R), \ t>0; \label{eq:1D-PDE-par} \\
& \bv{U}(x,0,y) = \bv{U}_0(x,y), \ y \in (y_L, y_R), \label{eq:1D-IC-par} \\
& \bv{U}(x_L,t,y) = \bv{U}_B^0(t,y), \ y \in (y_L, y_R), \label{eq:1D-BC0} \\
& \bv{U}(x_R,t,y) = \bv{U}_B^L(t,y), \ y \in (y_L, y_R). \label{eq:1D-BCL}
\end{align}

We can now treat the system \eqref{eq:1D-PDE-par}--\eqref{eq:1D-BCL} as a two-dimensional PDE in $(x,y)$ and discretize it following the principles of the multidimensional finite volume method.
Without loss of generality, we introduce the uniform grid in both physical and stochastic variables, with nodes $x_{i-1/2} = x_L + (i-1)\Delta x$ and
$y_{j-1/2} = y_L + (j-1)\Delta y$ that delimit the cells defined by $(i,j) = (x_{i-1/2},x_{i+1/2})\times(y_{j-1/2},y_{j+1/2})$.  We can then define the cell average for each cell $(i,j)$ by 
\begin{equation}
\bv{U}_{ij}(t) = \dfrac{1}{\Delta x |\Delta y|}\int\limits_{x_{i-1/2}}^{x_{i+1/2}}\int\limits_{y_{j-1/2}}^{y_{j+1/2}} \bv{U}(x,t,y)\mu(y)\,dxdy,
\end{equation}
where the volume of the cell in $y$-direction is defined as
\begin{equation}
|\Delta y| = \int\limits_{y_{j-1/2}}^{y_{j+1/2}} \mu(y)\,dxdy.
\end{equation}
Integrating the PDE \eqref{eq:1D-PDE-par} over the control volume $(x_{i-1/2},x_{i+1/2})\times(y_{j-1/2},y_{j+1/2})$ and applying integration by parts for the flux term results in the (still exact) formulation
\begin{multline}
 \Delta x |\Delta y| \diff{\bv{U}_{ij}}{t} + \left[ \int\limits_{y_{j-1/2}}^{y_{j+1/2}} \bv{F}\big(\bv{U}(x_{i+1/2},t,y)\big)\mu(y)\,dy - \int\limits_{y_{j-1/2}}^{y_{j+1/2}} \bv{F}\big(\bv{U}(x_{i-1/2},t,y)\big)\mu(y)\,dy \right] = \\ \int\limits_{x_{i-1/2}}^{x_{i+1/2}}\int\limits_{y_{j-1/2}}^{y_{j+1/2}} \bv{S}(\bv{U},y)\mu(y)\,dxdy.
\end{multline}

Thus applying the SFVM scheme to the PDE \eqref{eq:1D-PDE} with initial data \eqref{eq:1D-IC-par}, we obtain the following system of ODEs with respect to the cell averages:
\begin{align}
& \diff{\bv{U}_{ij}(t)}{t} + \dfrac{1}{\Delta x} \left[ \bv{F}_{i+1/2}^j(t) - \bv{F}_{i-1/2}^j(t) \right] = \bv{S}_{ij}(t), \label{eq:1D-ODE} \\
& \bv{U}_{ij}(0) = \dfrac{1}{\Delta x |\Delta y|}\int\limits_{x_{i-1/2}}^{x_{i+1/2}}\int\limits_{y_{j-1/2}}^{y_{j+1/2}} \bv{U}_0(x,y)\mu(y)\,dxdy. \label{eq:1D-ODE-IC}
\end{align}
We use the Gauss quadrature of appropriate order to compute the numerical fluxes in \eqref{eq:1D-ODE}:
\begin{equation}
\label{eq:quadrature}
\bv{F}_{i+1/2}^j(t) = \dfrac{1}{|\Delta y|}\int\limits_{y_{j-1/2}}^{y_{j+1/2}} \bv{F}\big(\bv{U}(x_{i+1/2},t,y)\big)\mu(y)\,dy \approx \sum\limits_{m=1}^M \hat{\bv{F}}\big(\tilde{\bv{U}}_{i+1/2,j}^{L,m}(t),\tilde{\bv{U}}_{i+1/2,j}^{R,m}(t)\big) w_m
\end{equation}
where $\hat{\bv{F}}(\cdot,\cdot)$ is any standard approximation of the flux (Godunov, Lax-Friedrichs, HLLC
flux, etc.).  A specific example of a numerical flux approximation is the Lax-Friedrichs flux,
\begin{equation} \label{eq:laxfriedrichsflux}
 \hat{\bv{F}}(\bv{U}_L,\bv{U}_R) = \dfrac{1}{2}(\bv{F}(\bv{U}_L) + \bv{F}(\bv{U}_R)) - \dfrac{\nu}{2}(\bv{U}_R-\bv{U}_L),
\end{equation}
where $\nu$ is a numerical viscosity coefficient. For our model of gas flow, the numerical viscosity coefficient is simply the wave speed $\nu = a$. The values $\tilde{\bv{U}}_{i+1/2,j}^{L,m}(t)$ and $\tilde{\bv{U}}_{i+1/2,j}^{R,m}(t)$ are obtained by reconstructing the solution inside the cells $(i,j)$ and $(i+1,j)$ using high order polynomials and evaluating these polynomials at the quadrature nodes $(x_{i+1/2},y_m)$. The source term in \eqref{eq:1D-ODE} is approximated by  
\begin{equation}
    \bv{S}_{ij}(t) = \dfrac{1}{\Delta x |\Delta y|}\int\limits_{x_{i-1/2}}^{x_{i+1/2}}\int\limits_{y_{j-1/2}}^{y_{j+1/2}} \bv{S}(\bv{U},y)\mu(y)\,dxdy.
\end{equation}

We refer the reader to detailed descriptions of ENO/WENO reconstruction \cite{shu99,zhang2016eno}. In the finite volume method, the solution is approximated by cell averages as a piecewise-constant function. Using cell averages as input data for the numerical fluxes results in at best a first order approximation. The ENO/WENO reconstruction enables the design of high-order finite volume methods in which the solution in each cell is represented as a polynomial reconstructed from cell averages, and this polynomial reconstruction serves as input data for the numerical fluxes. 

\subsection{Linear Total Variation Diminishing (TVD) Reconstruction} \label{sec:tvd_recon}

The following total variation diminishing (TVD) reconstruction is a simple example of a WENO scheme, which we present to demonstrate how to obtain cell interface data reconstructions $\tilde{\bv{U}}_{i+1/2,j}^{L,m}(t)$ and $\tilde{\bv{U}}_{i+1/2,j}^{R,m}(t)$ for the second-order SFV scheme. Given the cell averages $\bv{U}_{i-1,j}$, $\bv{U}_{ij}$, and $\bv{U}_{i+1,j}$, we can construct a linear (in $x$) function  ${\bv{U}}^L_j(x,t)$ inside the cell $(i,j)$ that satisfies the TVD (monotonicity) property. Similarly, we can construct a linear function ${\bv{U}}^R_j(x,t)$ that possesses TVD properties inside each cell $(i+1,j)$ by using the cell averages $\bv{U}_{ij}$, $\bv{U}_{i+1,j}$, and $\bv{U}_{i+2,j}$. Evaluating these functions at $x=x_{i+1/2}$ yields 
\begin{equation}
    \bv{U}^L_{i+1/2,j} = \bv{U}^L_j(x_{i+1/2},t), \quad \bv{U}^R_{i+1/2,j} = \bv{U}^R_j(x_{i+1/2},t)
\end{equation}
for each value of $j$.  Recall that the total variation of a function $u(x,t)$ is defined as
\begin{equation}
TV(u) = \int\Big|\pdiff{u}{x}\Big|\,dx, 
\end{equation}
and the total variation of the discrete function $u_h^n$ is
\begin{equation}
TV(u_h^n) = \sum\limits_j |u_{j+1}^n - u_j^n|.
\end{equation}
A discretization is said to be total variation diminishing (or monotonicity preserving) if it possesses the property
\begin{equation}
 TV(u_h^{n+1}) \leq TV(u_h^n).
\end{equation}

Consider the cell averages $\bv{U}_{i-1,j}$, $\bv{U}_{ij}$, $\bv{U}_{i+1,j}$. The goal is to construct a linear approximation of the solution satisfying the TVD property inside the cell $(i,j)$. The idea of the TVD reconstruction is to define candidate linear functions in the cell and then select the one with better properties. If the candidate function is defined as 
\begin{equation}
    \bv{U}_{ij}(x) = \bv{U}_{ij} + \sigma_{ij}(x-x_i),
\end{equation}
then the reconstruction reduces to finding an appropriate slope $\sigma_{ij}$. 
The first candidate reconstruction is obtained from $\bv{U}_{i-1,j}$ and $\bv{U}_{ij}$ and has slope $\sigma_{ij}^L = \frac{1}{\Delta x}(\bv{U}_{ij}-\bv{U}_{i-1,j})$, while the second candidate reconstruction uses $\bv{U}_{ij}$ and $\bv{U}_{i+1,j}$, and has slope $\sigma_{ij}^R = \frac{1}{\Delta x}(\bv{U}_{i+1,j}-\bv{U}_{ij})$. 

The TVD satisfying slope for the cell $(i,j)$ is then obtained by applying a slope limiter $\sigma_{ij} = \texttt{minmod}(\sigma_{ij}^L,\sigma_{ij}^R)$, where the \texttt{minmod} function is given by
\begin{equation}
\texttt{minmod}(a_1,a_2) = 
    \begin{cases}
        \rm{sign}(a_1)\min(|a_1|,|a_2|), \quad \text{if } \rm{sign}(a_1) = \rm{sign}(a_2), \\
        0, \quad \text{otherwise}.
    \end{cases}
\end{equation}

Next, we repeat the TVD linear reconstruction in the $y$-direction for each cell $(i,j)$ using the values $\bv{U}^L_{i+1/2,j-1}$, $\bv{U}^L_{i+1/2,j}$, and $\bv{U}^L_{i+1/2,j+1}$, as well as $\bv{U}^R_{i+1/2,j-1}$, $\bv{U}^R_{i+1/2,j}$, and $\bv{U}^R_{i+1/2,j+1}$. When we denote the resulting linear functions by $\hat{\bv{U}}^L_{i+1/2,j}(y,t)$ and $\hat{\bv{U}}^R_{i+1/2,j}(y,t)$, respectively,
Then we simply evaluate these linear functions to obtain the arguments for the numerical flux at a quadrature point $y_m$:
\begin{equation}
 \tilde{\bv{U}}_{i+1/2,j}^{L,m}(t) = \hat{\bv{U}}^L_{i+1/2,j}(y_m,t), \quad \tilde{\bv{U}}_{i+1/2,j}^{R,m}(t) = \hat{\bv{U}}^R_{i+1/2,j}(y_m,t).
\end{equation}

\begin{center}
\begin{tikzpicture}
\draw[gray, thick] (0,0) -- (5,0);
\draw[gray, thick] (0,3) -- (5,3);
\draw[gray, thick] (1,-1) -- (1,4);
\draw[gray, thick] (4,-1) -- (4,4);
\filldraw[black] (1,0) circle (2pt) node[anchor=north] {$x_{i-1/2}$};
\filldraw[black] (4,0) circle (2pt) node[anchor=north] {$x_{i+1/2}$};
\filldraw[black] (4,1) circle (2pt) node[anchor=north west] {$y_m$};
\draw[->] (2.5,1.5) -- (3.8,1);
\draw[->] (5.5,1.5) -- (4.2,1);
\draw (2.5,1.5) node[anchor=east] {$\tilde{\bv{U}}^L_{i+1/2}$};
\draw (5.5,1.5) node[anchor=west] {$\tilde{\bv{U}}^R_{i+1/2}$};
\end{tikzpicture}
\end{center}

The above reconstruction scheme can be generalized to the case when $y \in R^q$ for $q > 1$ in a straightforward manner by applying the multidimensional high-order ENO/WENO polynomial reconstruction procedure.

\section{Uncertainty Quantification for Pipeline Networks} \label{sec:networkuq}

We extend the SFV method to perform uncertainty quantification on a graph of PDEs, in order to apply the method to the setting of gas pipeline networks.  Such networks are referred to in previous studies as metric graphs \cite{facca2021fast} or flow networks \cite{misra2020monotonicity}.  We suppose the setting of the underlying deterministic system to be a graph consisting of a set of vertices $\cV$ that are connected by a set $\cE$ of oriented edges. The incoming and outgoing neighborhoods of a node $j$ are denoted by $\partial_{+}j$ and $\partial_{-}j$, respectively, and are defined as
\begin{align}
\partial_{+}j&=\left\{ i\in \cV\mid(i,j)\in \cE\right\} \\
\partial_{-}j&=\left\{ k\in \cV\mid(j,k)\in \cE\right\}.
\end{align}
Every edge $(i,j)\in \cE$ is associated with a spatial dimension on the interval $I_{ij}=[0,L_{ij}]$, where $L_{ij}$ is the edge length.  The flow dynamics on each edge $(i,j)\in\cE$ are described by per area mass flux $q_{ij}(t,x)$ and density $\rho_{ij}(t,x)$.  Each node $j$ has a unique nodal density $\rho_{j}(t)$, and flow leaving the network at that node is denoted $d_{j}(t)$.

We suppose that the density and flux dynamics on an edge $(i,j)$ evolve according to
\begin{subequations}\label{eq:gas1}
\begin{align}
& \pdt \rho_{ij} + \pdx q_{ij} = 0, \\
& \pdt q_{ij} + a^2\pdx \rho_{ij} = -\dfrac{\lambda_{ij}}{2D_{ij}}\dfrac{q_{ij}|q_{ij}|}{\rho_{ij}},
\end{align}
\end{subequations}
where each edge has parameters $\lambda_{ij}$ and $D_{ij}$. Here we use ideal gas modeling, and suppose that the wave speed $a$ is constant throughout the network.  The SFV methodology could be extended to non-ideal gas flow using local nonlinear transformations to adjust this parameter. We use a shorthand notation for values of the state variables at the edge boundaries:
\begin{align}
\underline{\rho}_{ij}(t)\triangleq\rho_{ij}(t,0), \quad \overline{\rho}_{ij}(t)\triangleq\rho_{ij}(t,L_{ij}), \label{eq:end_p_def} \\
\underline{q}_{ij}(t)\triangleq q_{ij}(t,0), \quad \overline{q}_{ij}(t)\triangleq q_{ij}(t,L_{ij}). \label{eq:end_q_def}
\end{align}
At each vertex $j$, the flow and density values at the endpoints of adjoining edges must satisfy certain compatibility conditions.  First, a Kirchhoff-Neumann property of flow conservation is ensured through nodal continuity equations
\begin{align}
\sum_{i\in\partial_{+}j}X_{ij}\overline{q}_{ij}(t)- \sum_{k\in\partial_{-}j}X_{jk}\underline{q}_{jk}(t)= d_j(t), \quad \forall j\in\cV, \label{eq:in_nodal_continuity} 
\end{align}
where $X_{ij}$ denotes the cross-sectional area of pipe $(i,j)$.  In addition, we incorporate modeling of gas compressors, which are a set of elements $\mathcal C$ that are used to actuate flow throughout the network.  Each compressor $c\in\mathcal C$ is located at a node $j$ and affects an increase between the nodal density $\rho_j(t)$ and the boundary density $\underline{\rho}_{jk}(t)$ where $(j,k)\in\mathcal E$ is a pipe oriented outward from node $j$.  The action of such a compressor is described by the multiplicative ratio $\alpha_{jk}(t)$, as
\begin{align}
\underline{\rho}_{jk}(t) = \alpha_{ij}(t)\rho_{j}(t), \label{eq:in_pressure_comp} 
\end{align}
We suppose that instantaneous state of the system at time $t=0$ is specified by initial density and mass flux profiles
\begin{align}
\!\!\! \rho_{ij}(0,x)=\rho_{ij}^{0}(x), \,\, q_{ij}(0,x)=q_{ij}^{0}(x), \quad \forall (i,j)\in\cE. \label{eq:in_initial_condition}
\end{align}

\subsection{Conditions at Junctions}

Here we consider in detail the numerical treatment of uncertainty propagation through nodes of flow networks.  The conditions that must be enforced include flow conservation \eqref{eq:in_nodal_continuity} and density compatibility \eqref{eq:in_pressure_comp}.  We extend the above modeling to account for the conservation of the distributions of conserved quantities.  Consider a localized IBVP for node $j\in\mathcal V$ and involving all edges $e\in \mathcal E_j \equiv \partial_-j \cup \partial_+ j \subset \mathcal E$, which can be stated as
\begin{subequations}
\begin{align}
&\pdiff{\bv{U}_{ij}}{t} + \nabla_x\cdot\bv{F}(\bv{U}_{ij}) = \bv{S}(\bv{U}_{ij}), \,\,\, x_{ij} \in [0,L_{ij}],  t > 0, \,\,\, & \forall e\in\mathcal E_j, \\
&\bv{U}_{ij}(0,x) = \bv{U}^0_{ij}(x_{ij}), & \forall e\in\mathcal E_j.
\end{align} \label{eq:PDEs-vertex}
\end{subequations}
The material and geometric properties (i.e., $\lambda_{ij}$ and $D_{ij}$) can be different for each edge $e\in\mathcal E_j$. 
The set of solutions $\bv{U}_{ij}(t,x_{ij})$ has to satisfy the coupling conditions at the common node $j$.  Without loss of generality, we may assume that all edges are oriented outwards from $j$, such that $\partial_+j=\emptyset$, and where the incoming nodes are indexed $k=1,\ldots,K$.  The coupling conditions for node $j$ can then be stated as
\begin{equation}
\label{Eq:coupling-V}
\bv{\Phi}(\bv{U}_{j1}(t,0),\dots,\bv{U}_{jK}(t,0)) = \bv{0}, \ t > 0.
\end{equation}
Equations \eqref{eq:PDEs-vertex}--\eqref{Eq:coupling-V} define a Junction Riemann problem that extends a classical Riemann problem corresponding to $N = 2$.  

For a physical flow network such as a gas pipeline system, the conditions at the node $j$ will consist in continuity of pressure and conservation of flow, so that in the above case where $\partial_+j=\emptyset$ and incoming nodes are indexed $k=1,\ldots,K$, the coupling conditions are 
\begin{equation}
\label{Eq:coupling-gas}
\bv{\Phi} = 
    \begin{pmatrix}
        \rho_{jk} - \alpha_{jk} \rho_j, \ \forall k \in \delta_-j  \\
        \sum_{k \in \delta_-j} X_{jk}q_{jk} (t,0) -d_j(t)
    \end{pmatrix}
= \bv{0}.
\end{equation}

\subsection{Classical Riemann problem}

The SFV method utilizes a cell-centered approximation of density and flow. These cell-centered values are updated using the values of numerical fluxes at the cell interfaces. A node corresponding to a given junction coincides with the cell interfaces of the first or last cells in the discretization of each pipe. Therefore, the numerical flux through those cell interfaces must take into account the physical compatibility conditions. The numerical fluxes in the finite volume method are typically computed by considering Riemann problems at the interface between the two cells. A classical Riemann problem is a one-dimensional Cauchy problem with piecewise-constant initial data:
\begin{subequations}
\begin{align}
& \pdiff{\bv{U}}{t} + \nabla_x\cdot\bv{F}(\bv{U}) = \bv{0}, \ t > 0; \\
& \bv{U}(x,0) = 
    \begin{cases}
     \bv{U}_L, \quad \text{if} \quad x < 0, \\
     \bv{U}_R, \quad \text{if} \quad x > 0.
    \end{cases}
\end{align} \label{RP-classical}
\end{subequations}

\begin{figure}[h!]
    \centering
    \includegraphics[scale=0.5]{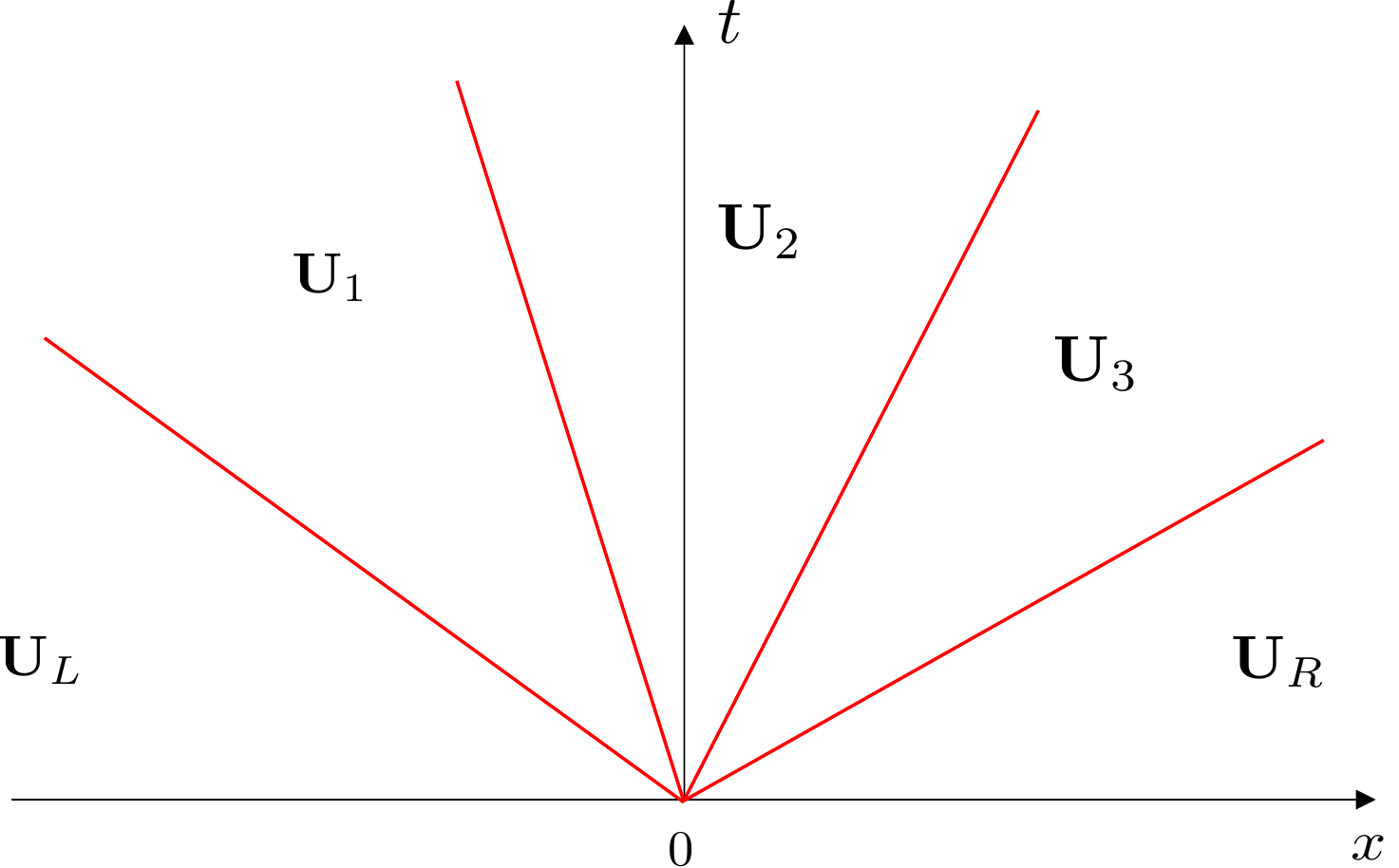}
    \caption{Structure of the solution of the Riemann problem on $(x,t)$ plane}
    \label{Fig:RP-waves}
\end{figure}

The solution of the Riemann problem \eqref{RP-classical} is a self-similar function $\bv{U} = \bv{D}(x/t)$. For a hyperbolic system of $m$ equations, the solution on the $(x,t)$ plane consists of $m+1$ regions of constant values of $\bv{U}$, which are separated by $m$ characteristic waves, as illustrated in Fig.~\ref{Fig:RP-waves} for $m=4$. 
The physical variable values across each of the $m$ waves are connected by means of generalized Riemann invariants $I_l$, so that $I_l(\bv{U}) = \const$, $l = 1,\dots,m$.  There exist a number of analytical and numerical methods to obtain the solution of \eqref{RP-classical}. For the purpose of determining the numerical fluxes in finite volume and related methods, we are only interested in the solution along the $t$-axis, that is, $\bv{D}(0)$.  In the following section, we describe a generalization of the classical Riemann problem at junctions.

\subsection{Junction Riemann problem (JRP)}

\begin{figure}[h!]
    \centering
    \includegraphics[scale=0.5]{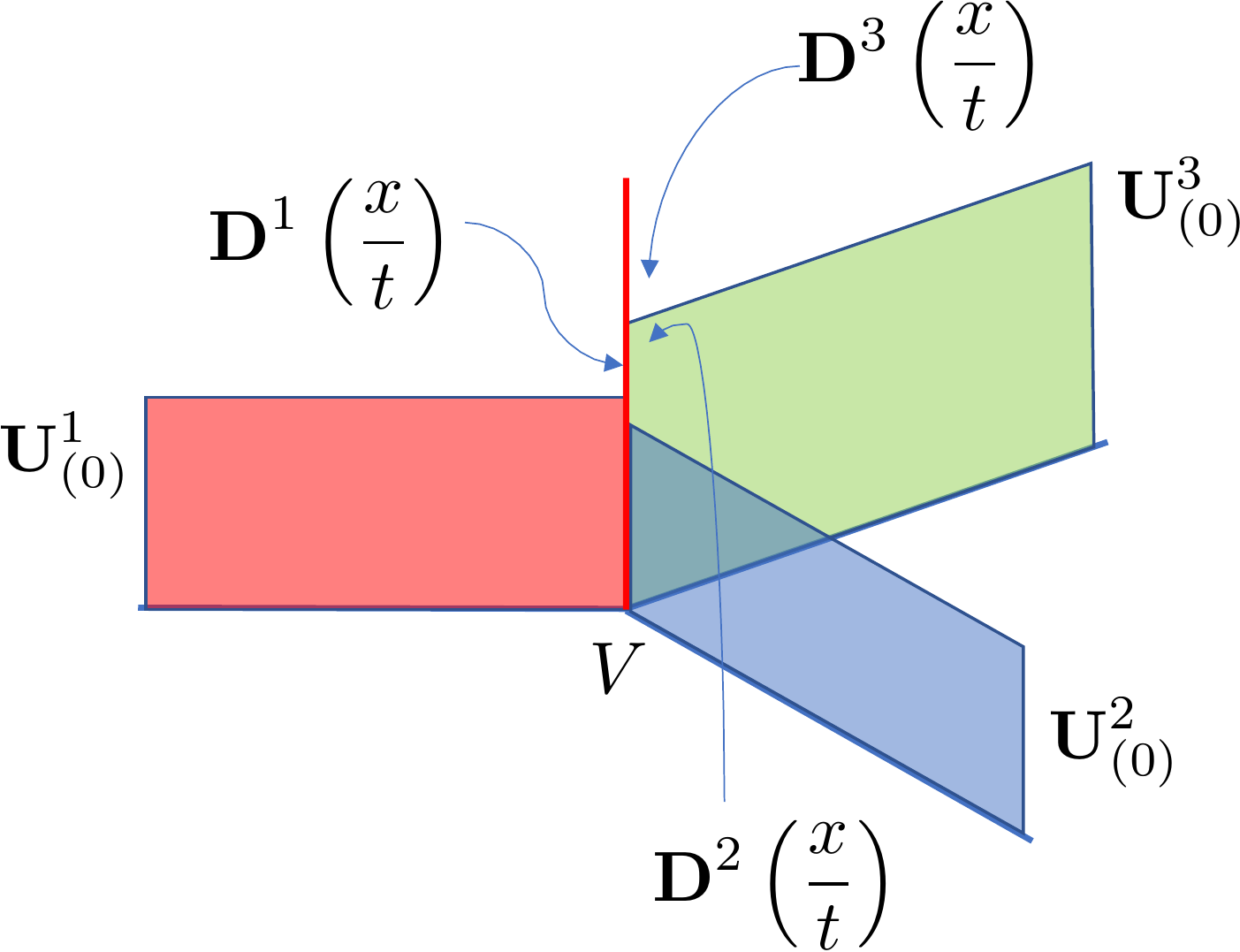}
    \caption{Junction Riemann problem for $3$ pipes}
    \label{Fig:J-CRP}
\end{figure}

A generalization of the Riemann problem for network junctions is illustrated in Fig.~\ref{Fig:J-CRP}, where three pipes are connected at a vertex $V$, see also \cite{JRP2016,UQBF2019}.
The numerical discretization of the gas flow parameters that we apply on each pipe is based on the finite volume approach, and is therefore cell-centered. The last cell of the incoming pipe or the first cell of the outgoing pipe will contain an interface that coincides with the vertex $V$.

The one dimensional gas flow dynamics \eqref{Eq:gas0} is a linear hyperbolic system of two equations, and thus it has two characteristic waves. The corresponding Riemann invariants are 
\begin{align}
& I_1(\rho,q) = \rho + \frac{1}{a}q = \const, \\
& I_2(\rho,q) = \rho - \frac{1}{a}q = \const.
\end{align}
For the collection of pipes adjoining a vertex $V$, let us denote the solution of the JRP in the $k$-th pipe along the $x$-axis by $\rho^k_*$ and $q^k_*$. For the three pipe junction example illustrated in Fig.~\ref{Fig:J-CRP}, the solutions of the JRP is obtained by solving the following system, which combines the Riemann invariants and constraints:
\begin{subequations}
\begin{align}
\label{Eq:JRP-system}
& q^1_* = q^2_* + q^3_*, \\
& a^2 \rho^1_* = a^2 \rho^2_*, \\
& a^2 \rho^1_* = a^2 \rho^3_*, \\
& \rho^1_* - \rho^1_{(0)} + \frac{1}{a} (q^1_* - q^1_{(0)}) = 0, \\
& \rho^2_* - \rho^2_{(0)} - \frac{1}{a} (q^2_* - q^2_{(0)}) = 0, \\
& \rho^3_* - \rho^3_{(0)} - \frac{1}{a} (q^3_* - q^3_{(0)}) = 0.
\end{align}
\end{subequations}
This is a linear system of $6$ equations with respect to $6$ unknowns. Having determined the pairs $U^k_* = (\rho^k_*,q^k_*)$, we can calculate numerical fluxes at the vertex $V$ using, e.g., the Lax-Friedrichs or Godunov flux $F_V^k = F(U^k_*)$.

\section{Numerical Examples} \label{sec:computation}

Here we apply the uncertainty quantification method described above first to flow through a single pipe, and then to flow through a small test network. We suppose that density is fixed at the left (inlet) end of the pipe, and that the source of uncertainty is in the mass flow boundary condition at the right (outflow) end of the pipe. 

\subsection{Uncertainty Quantification for Flow in One Pipe} \label{sec:onepipemodel}

We consider the gas flow PDEs given in Equation \eqref{Eq:gas0}, with steady-state initial conditions given by
\begin{subequations} \label{eq:onepipeic}
\begin{align}
\label{eq:ic-rho}
& \rho(0,x) = \sqrt{s(0)^2 - \dfrac{\lambda}{a^2 D} q_0 |q_0|x}, \\
\label{eq:ic-q}
& q(0,x) = q_0,
\end{align}
\end{subequations}
and time-varying boundary conditions given for $t\geq 0$ by
\begin{subequations} \label{eq:onepipebc}
\begin{align}
& \rho(t,0) = s(t), \\
& q(t,L) = d(t).
\end{align}
\end{subequations}
The parameters in the PDE system and the initial and boundary conditions are those used in two previous studies on deterministic methods \cite{herty2010new,zlotnik2015model}.  The parameters are $\lambda = 0.011$, $a = 377.9683$ m/s, $D = 0.5$ m, $L = 100 \times 10^3$ m, $T = 3600 \times 12 s$; the nominal values of flow and density are $q_0 = 289$ kg/m$^2$/s and $\rho_0 = 45.4990786148$ kg/m$^3$, respectively; and the boundary conditions for the underlying deterministic problem are
\begin{subequations} \label{eq:onepipebcvals}
\begin{align}
& s(t) = \rho_0 \left(1+\frac{1}{10}\sin(6\pi t/T)\right), \\
& d(t) = q_0 \left(1 + \frac{1}{10}\sin(4\pi t/T)\right). \label{eq:onepipebcvals_q}
\end{align}
\end{subequations}
We extend the above problem formulation to a stochastic setting by adding a source of model uncertainty to the flow boundary condition. In particular, we suppose that the right boundary condition depends on the random variable $\omega \in \Omega$, so that mathematically we can write $q(L,t) = \hat{d}(t;\omega)$. We consider several scenarios for $\hat{d}(t;\omega)$ to examine solutions for both interval and inter-temporal uncertainty, as illustrated conceptually in Figure \ref{fig:uncertainties-peak}.

\subsubsection{Interval uncertainty in the boundary condition} \label{sec:example1}

\begin{figure}[h!]
\includegraphics[width=0.5\textwidth]{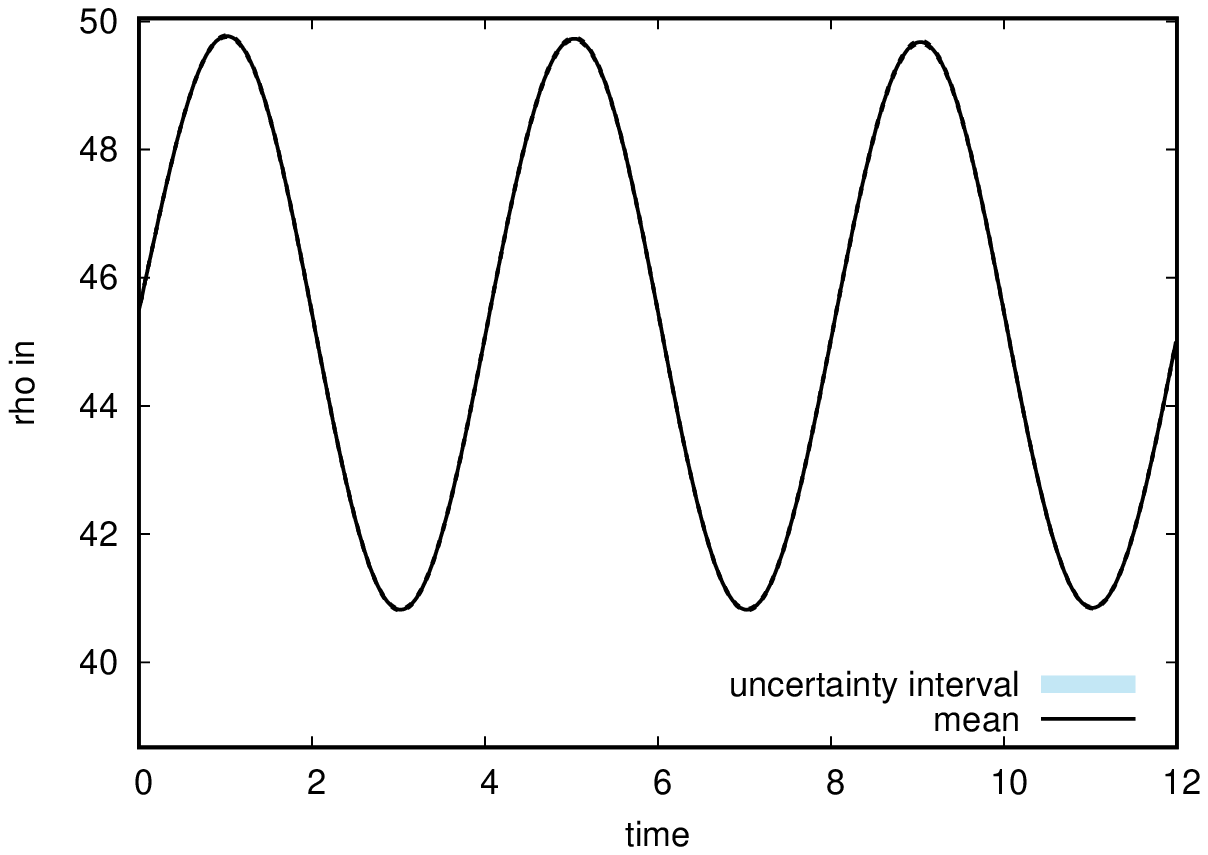}   \includegraphics[width=0.5\textwidth]{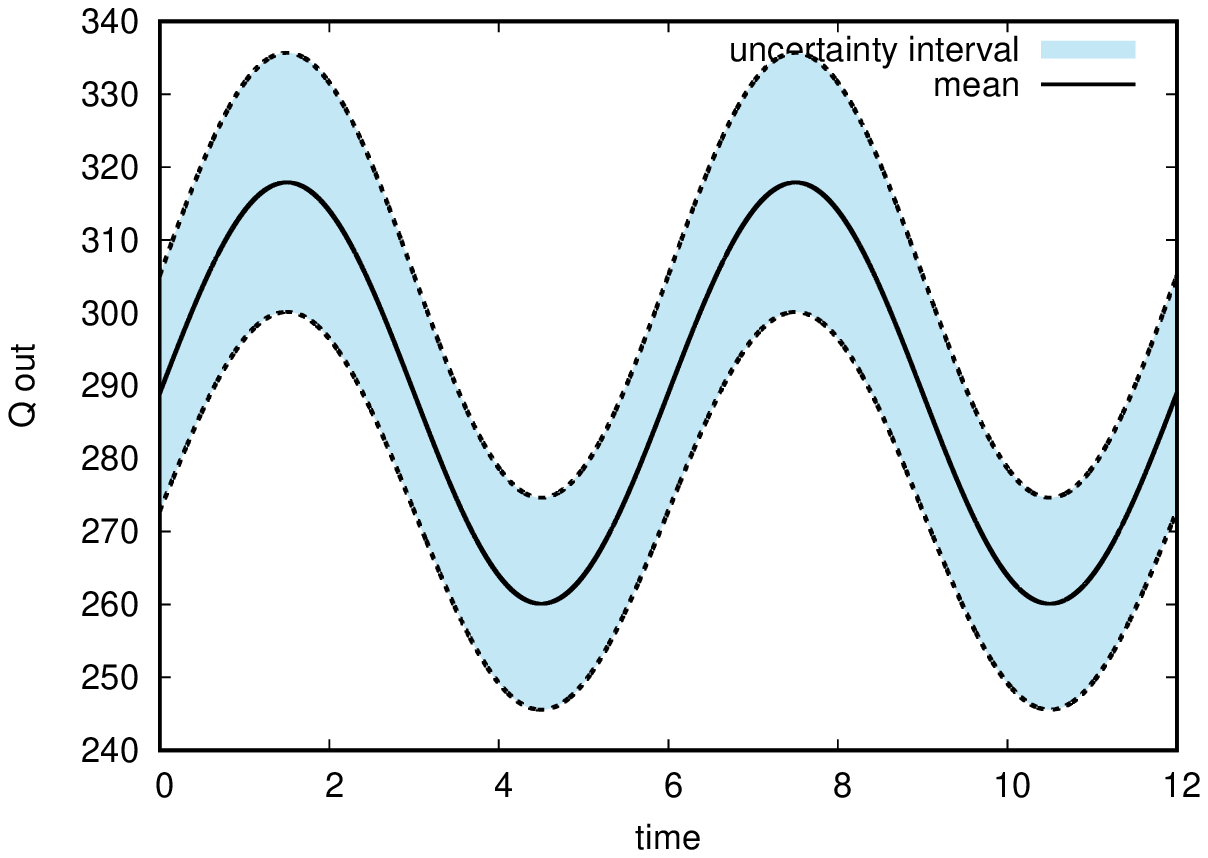} 
\caption{Boundary conditions used as inputs for simplest one pipe simulation. Left: density $\rho(t,0)$ at the pipe inlet; Right: outflow $q(t,L)$ from the pipe. The shaded region is one standard deviation above and below the mean.}
\label{fig:comp_onepipe_1_bc}
\end{figure}

\begin{figure}[h!]
\includegraphics[width=0.5\textwidth]{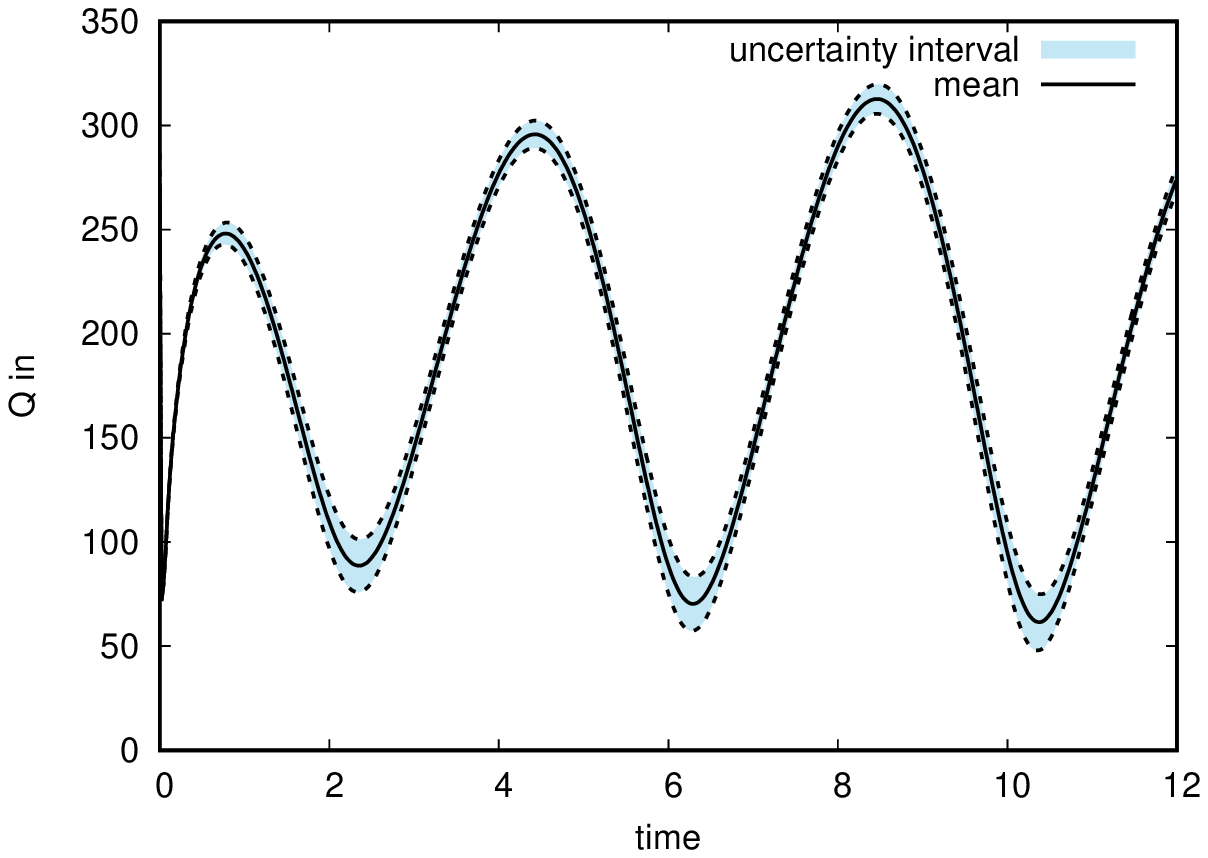}   \includegraphics[width=0.5\textwidth]{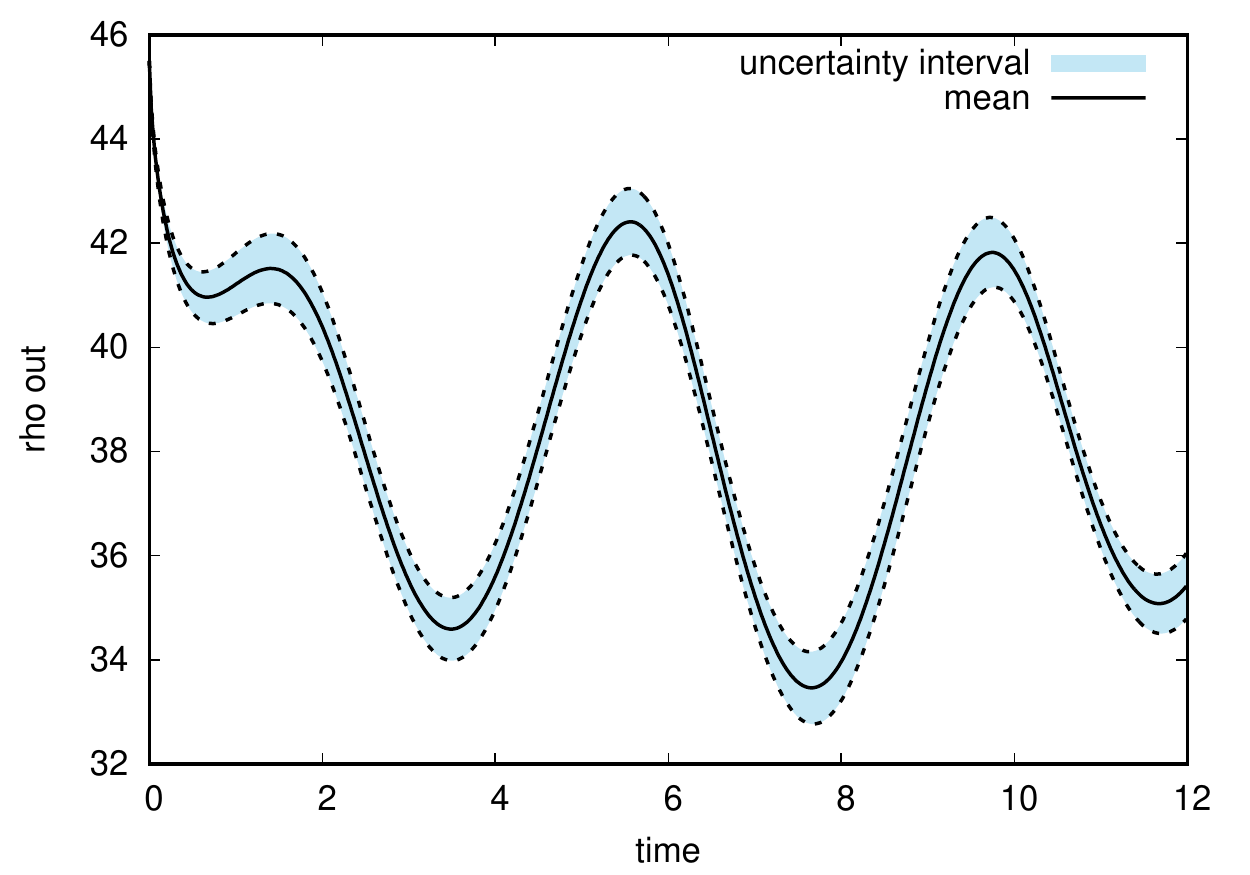} 
\caption{Uncertainty quantification results for simplest one pipe simulation. Left: inflow $q(t,0)$ at pipe inlet; Right: discharge density $\rho(t,L)$ at pipe outlet;  The shaded region is one standard deviation above and below the mean.}
\label{fig:comp_onepipe_1_output}
\end{figure}

Recall the notion of interval uncertainty as illustrated in Figure \ref{fig:uncertainties-peak}, which considers the possibility of time-varying parameters that lie between upper and lower bounding functions.  We first consider the simplest scenario, in which the nominal flow parameter $q_0$ in Equation \eqref{eq:onepipebcvals_q} is randomly defined as $q_0 \sim 289Y(\omega)$, where $Y(\omega) \sim U[0.9,1.1]$ is uniformly distributed.  The boundary conditions are shown in Figure \ref{fig:comp_onepipe_1_bc}, where the density at the inlet is deterministic and the interval of uncertainty for the outlet flow is indicated.  The resulting solutions at the boundaries are then shown in Figure \ref{fig:comp_onepipe_1_output}, and the uncertainty shown in the inlet flow and outlet density is one standard deviation above and below the mean.


\subsubsection{Interval uncertainty with a normal distribution} \label{sec:example2}

\begin{figure}[h!]
\includegraphics[width=0.5\textwidth]{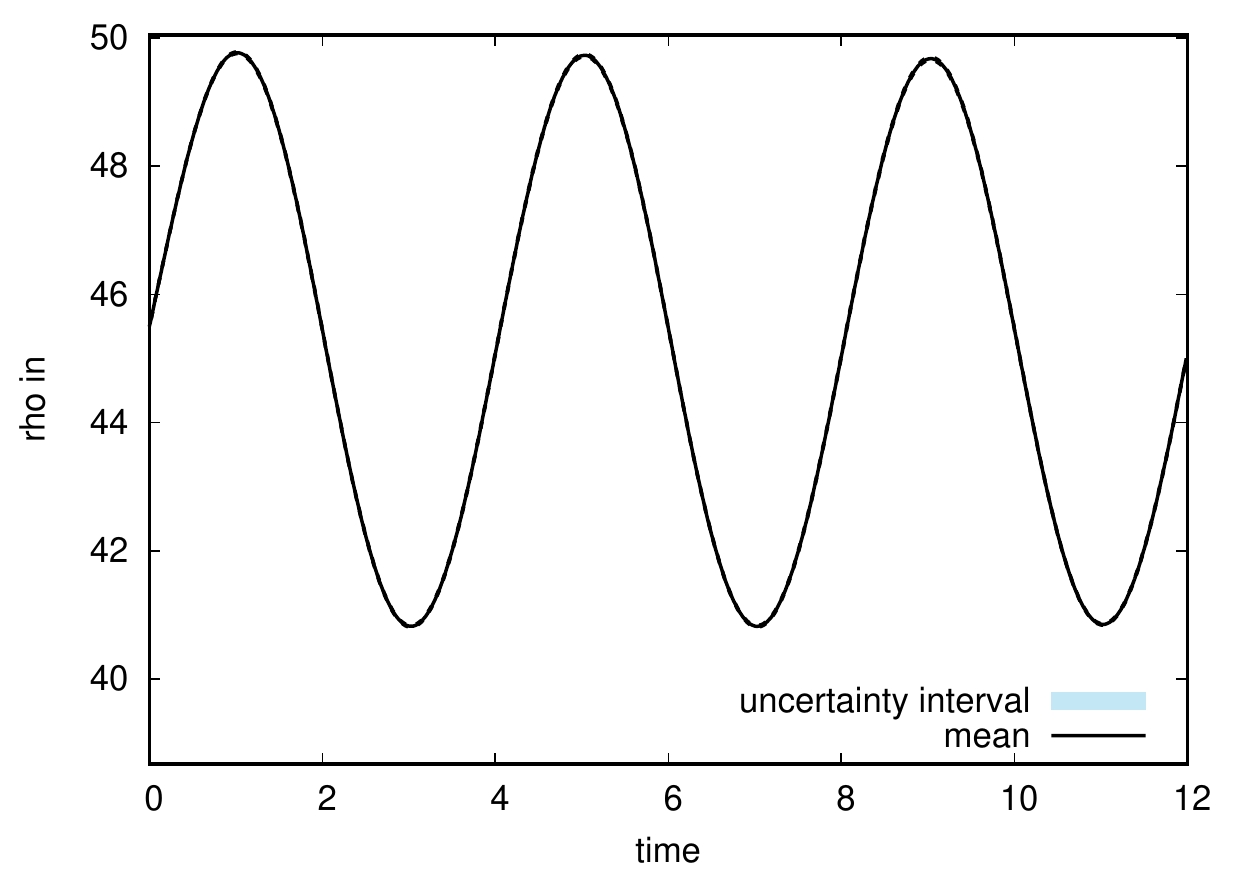}   \includegraphics[width=0.5\textwidth]{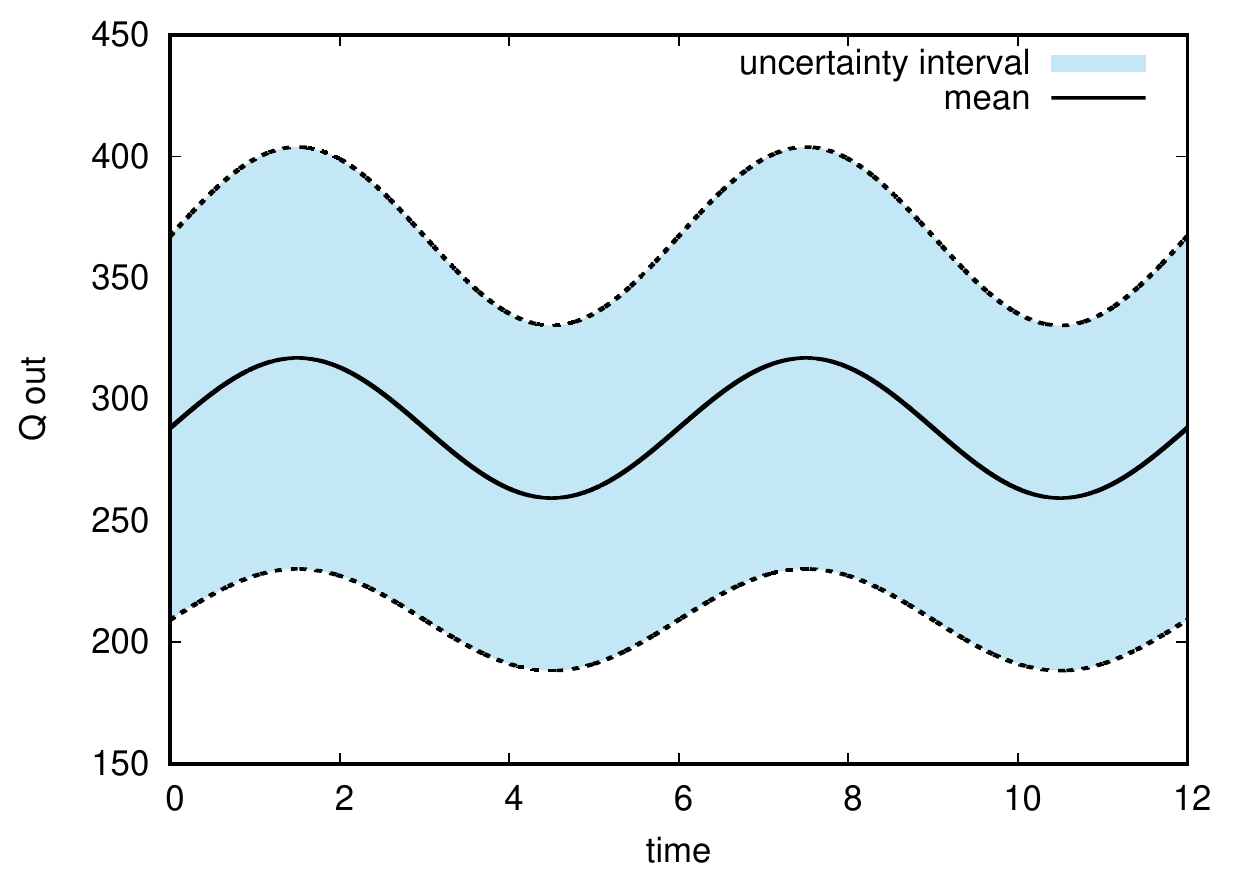} 
\caption{Boundary conditions used as inputs for one pipe simulation with Gaussian uncertainty. Left: density $\rho(t,0)$ at the pipe inlet; Right: outflow $q(t,L)$ from the pipe. The shaded region is one standard deviation above and below the mean.}
\label{fig:comp_onepipe_2_bc}
\end{figure}

\begin{figure}[h!]
\includegraphics[width=0.5\textwidth]{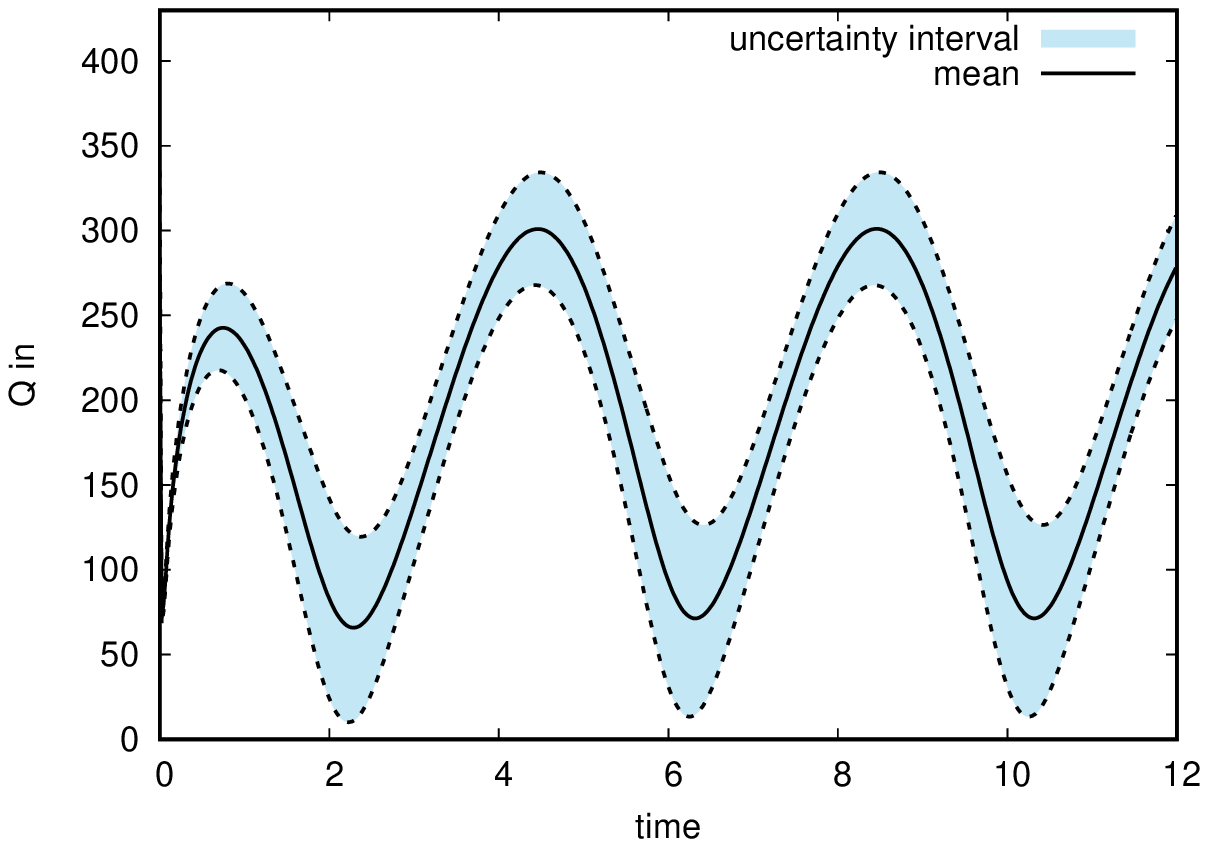}   \includegraphics[width=0.5\textwidth]{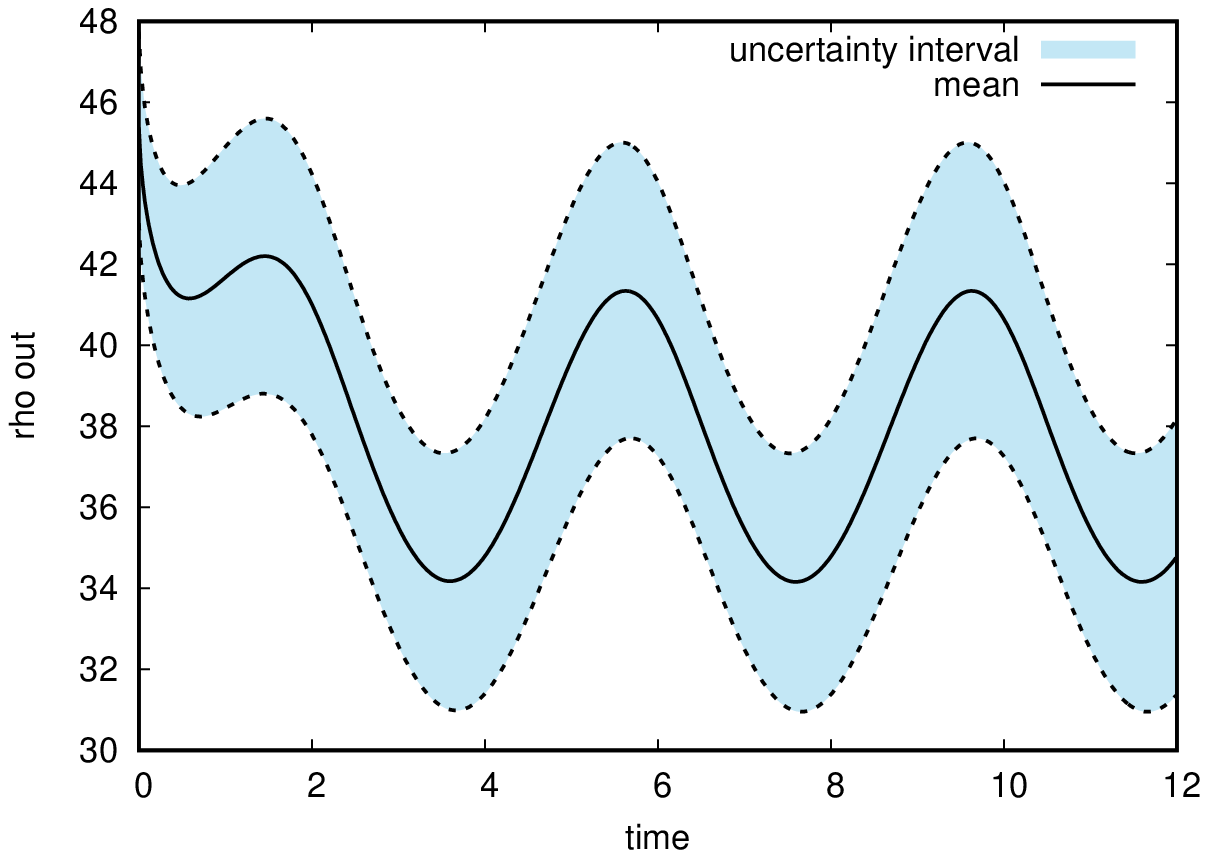} 
\caption{Uncertainty quantification results for one pipe simulation with Gaussian uncertainty. Left: inflow $q(t,0)$ at pipe inlet; Right: discharge density $\rho(t,L)$ at pipe outlet;  The shaded region is one standard deviation above and below the mean.}
\label{fig:comp_onepipe_2_output}
\end{figure}

Next, we consider the example as described in Section \ref{sec:example1}, but instead the nominal flow parameter $q_0$ in Equation \eqref{eq:onepipebcvals_q} is $q_0 \sim 289 Y(\omega)$, where  $Y(\omega) \sim N(0.25,1)$ is a normally distributed random variable.  The boundary conditions are shown in Figure \ref{fig:comp_onepipe_2_bc}, where the density at the inlet is deterministic and the interval of uncertainty for the outlet flow is indicated.  The resulting solutions at the boundaries are then shown in Figure \ref{fig:comp_onepipe_2_output}, and the uncertainty shown in inlet flow and outlet density is one standard deviation above and below the mean.

\subsubsection{Convergence analysis of the SFV method} \label{sec:convergence}

\begin{figure}[h!]
\includegraphics[width=0.5\textwidth]{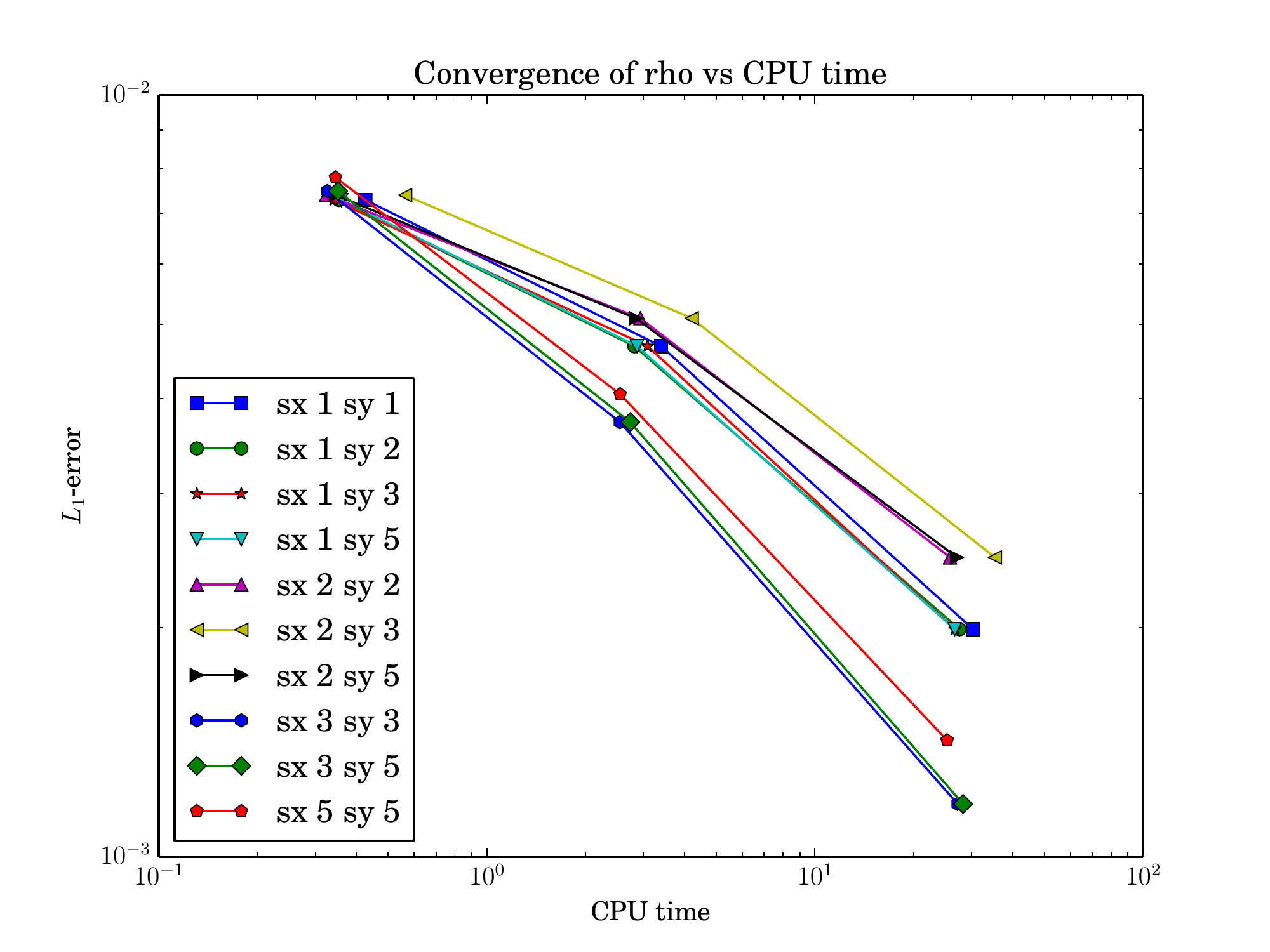}   \includegraphics[width=0.5\textwidth]{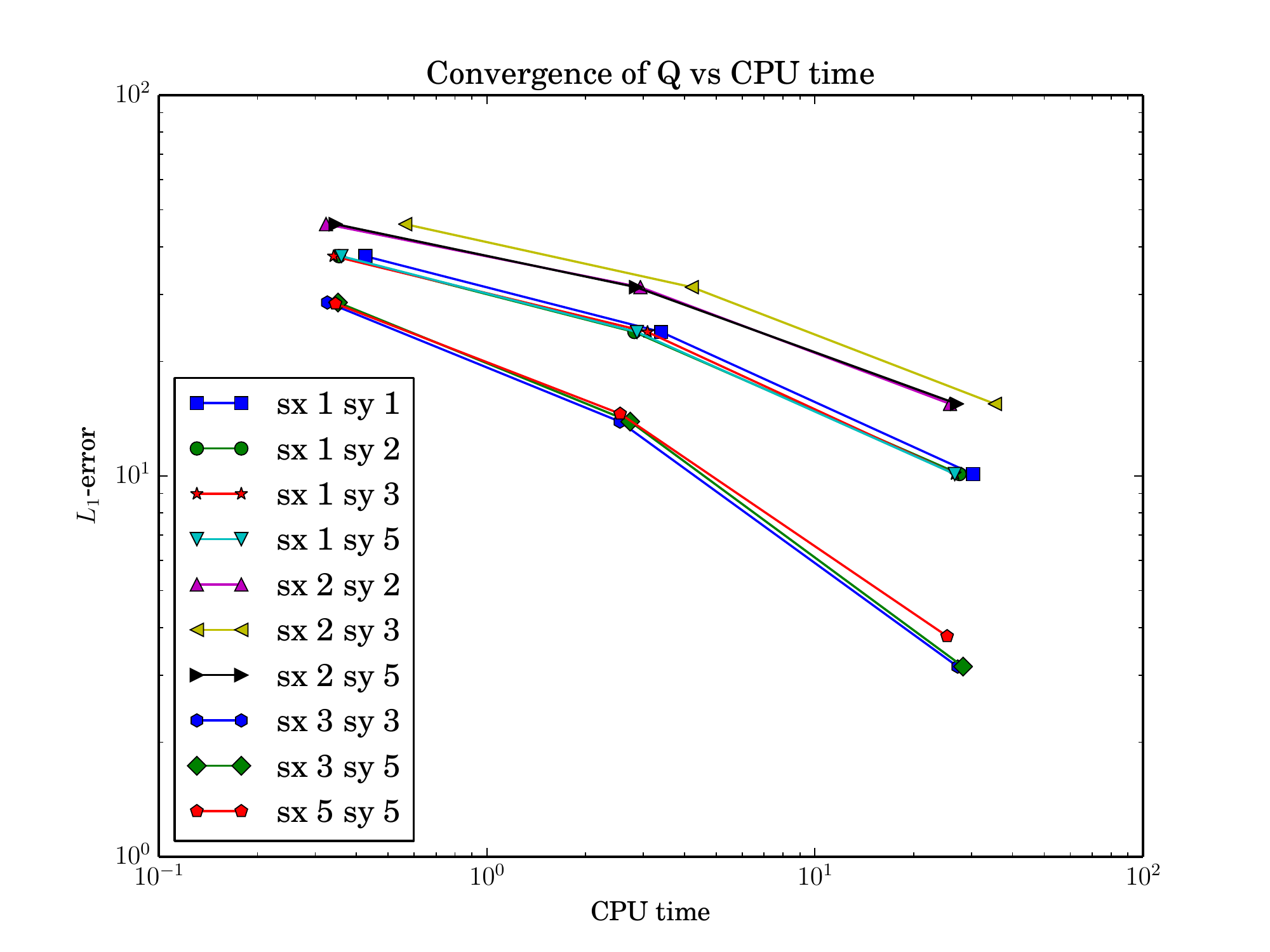} 
\caption{Computational and computation time results.  The convergence of the SFV computation is shown for outlet density $\rho$ (left) and inlet flow $q$ (right).  The plots show $L_1$ error vs. CPU time.}
\label{fig:onepipe_convergence}
\end{figure}

For the convergence study, we solve equations \eqref{Eq:gas0} with initial and boundary conditions given by \eqref{eq:ic-rho}--\eqref{eq:onepipebc} and random $q_0 \sim 289 Y(\omega)$, where  $Y(\omega) \sim U[0.9,1.1]$ is a uniformly distributed random variable. We repeat the computation on a series of meshes with $N_x = [4,8,16,32]$ and corresponding $N_y = [1,2,4,8]$, and varying orders of approximation in both physical and stochastic space. Figure~\ref{fig:onepipe_convergence} shows the error in $L_1$-norm as a function of the CPU time required to complete the simulation on each of the meshes. The solid lines correspond to different combinations of reconstruction in physical and stochastic coordinates: for example, "sx 3 sy 5" in the legend means that WENO scheme of order $3$ was used in the physical space and WENO scheme of order $5$ --- in the stochastic space. The plots in Figure~\ref{fig:onepipe_convergence} clearly show that the computational efficiency of the method increases with the order of reconstruction (i.e. the CPU time needed to achieve the desired accuracy decreases).

\subsubsection{Inter-temporal uncertainty in the boundary condition} \label{sec:example3}

\begin{figure}[h!]
    \centering
    \includegraphics[width=0.45\textwidth]{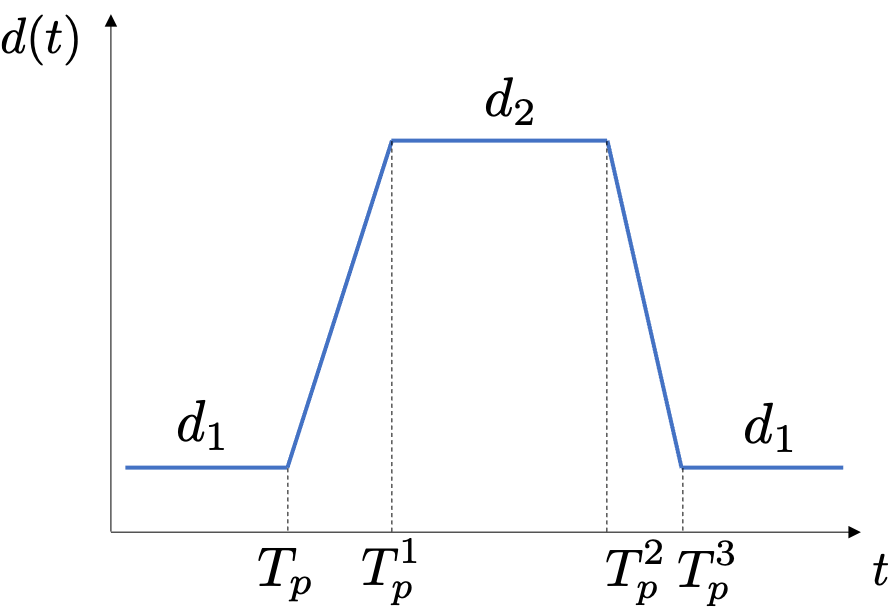}
    \caption{Example withdrawal rate $d(t)$ featuring a temporary perturbation.}
    \label{fig:perturb_dt}
\end{figure}

The results in the above examples with interval uncertainty can be obtained in a straightforward manner using other methods.  The key advantage of the SFV approach in our study is the ability to quantify the effect of inter-temporal uncertainty, as shown in the next example.  Recall that, as illustrated in Figure \ref{fig:uncertainties-peak}, we define inter-temporal uncertainty as a perturbative increase in flow with a specified duration and an uncertain initial time.  In the simplest example of inter-temporal uncertainty quantification for gas flow, we suppose that the time $T_p$ when the perturbation begins is a uniformly distributed random variable $Y(\omega)$, and a perturbation to the boundary flow function $d(t)$ is modeled using the following function:
\begin{equation}
d(t) = \left\{
    \begin{array}{ll}
     d_1(t), & \text{if} \quad t < T_p, \\
     (d_2-d_1)/(T_p^1 - T_p)(t-T_p) + d_1(t), & \text{if} \quad t \in (T_p,T_p^1) \\
     d_2, & \text{if} \quad t \in (T_p^1,T_p^2) \\
     (d_1-d_2)/(T_p^3 - T_p^2)(t-T_p^2) + d_2, & \text{if} \quad t \in (T_p^2,T_p^3) \\
     d_1(t), & \text{if} \quad t > T_p^3.
    \end{array}\right.
\end{equation}
An example boundary flow function $d(t)$ constructed using constant functions $d_1(t)\equiv d_1$ and $d_2(t)\equiv d_2$ is illustrated in Fig.~\ref{fig:perturb_dt}.

In the computations below, we suppose that $d_1(t)$ is given by the nominal deterministic boundary condition in Equation \eqref{eq:onepipebcvals_q}.  We then suppose that $d_2(t) = 3\cdot d_1(t)$, and the perturbation timing is specified by the initial time $T_p(\omega) = 3600\cdot(1 + Y(\omega))$ seconds, and the duration is $dT_p = 5\cdot 3600$ seconds, so that $T_p^1 = T_p + 0.1\cdot dT_p$, $T_p^2 = T_p + dT_p - 0.1\cdot dT_p$, and $T_p^3 = T_p + dT_p$. For the first computation, we assume that $Y(\omega) \sim U[0,2]$, and therefore $T_p$ has a compact support. 
\begin{figure}[h!]
\includegraphics[width=0.5\textwidth]{Fig/rho_vs_time_in_t12h_q0}   \includegraphics[width=0.5\textwidth]{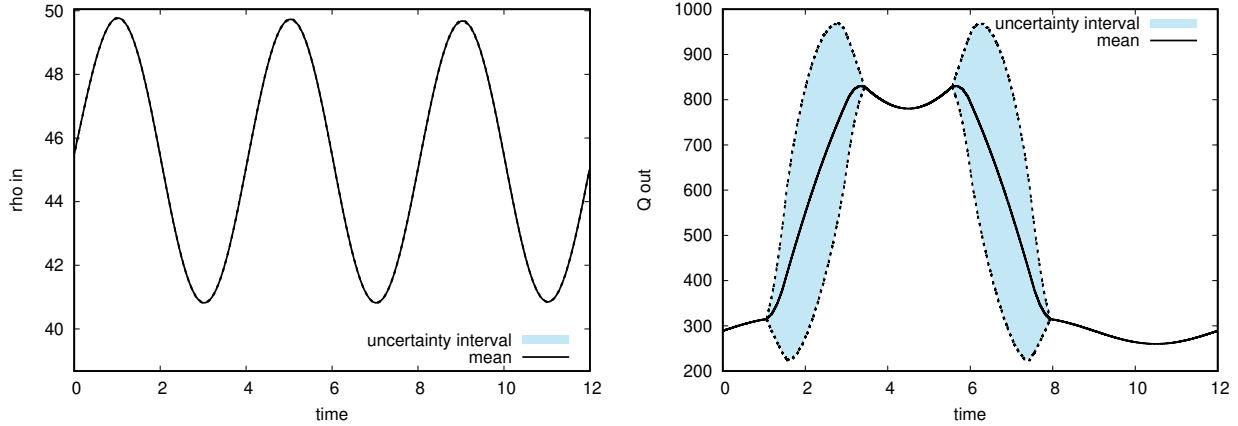} 
\caption{Boundary conditions used as inputs for one pipe simulation with limited inter-temporal uncertainty.  Left: inlet density $\rho(t,0)$; Right: outlet flow $q(t,L)$.  The shaded region is as previously defined.}
\label{fig:comp_onepipe_3_bc}
\end{figure}
\begin{figure}[h!]
\includegraphics[width=0.5\textwidth]{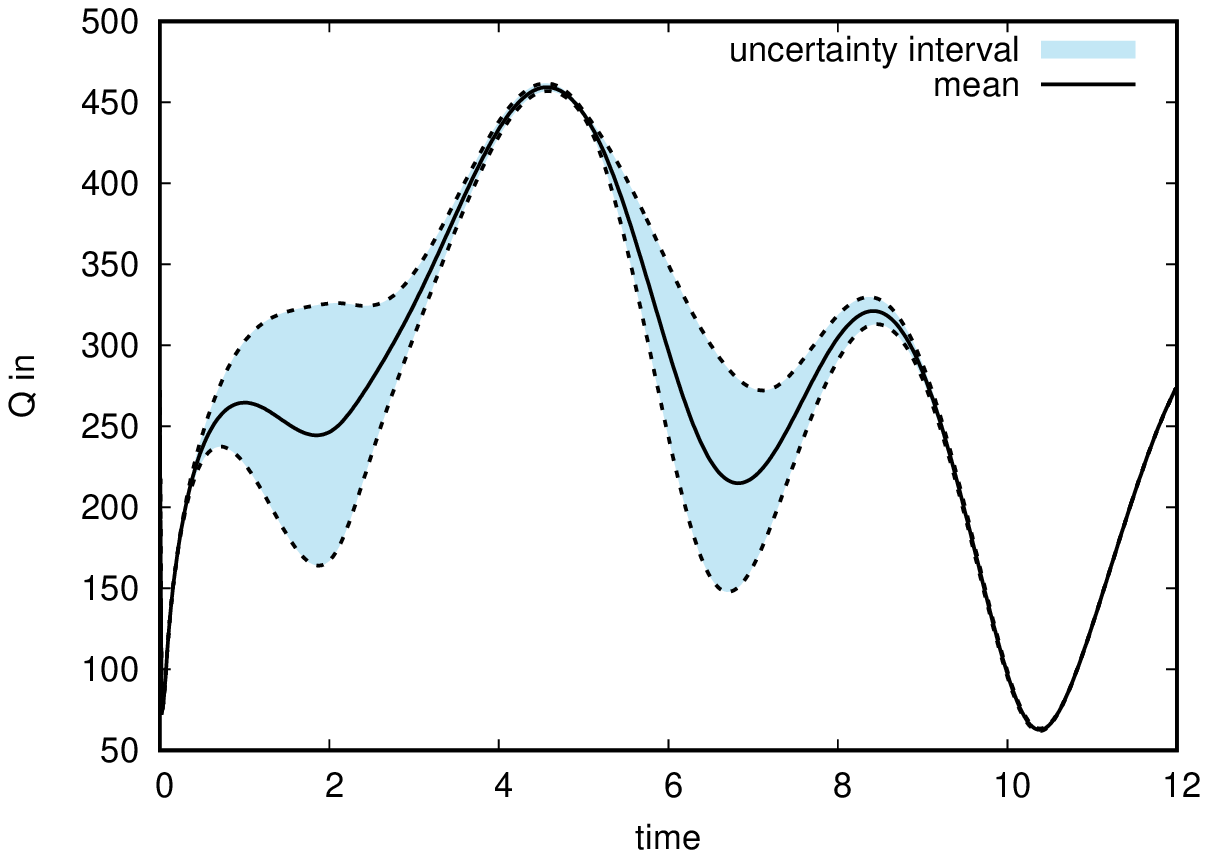}   \includegraphics[width=0.5\textwidth]{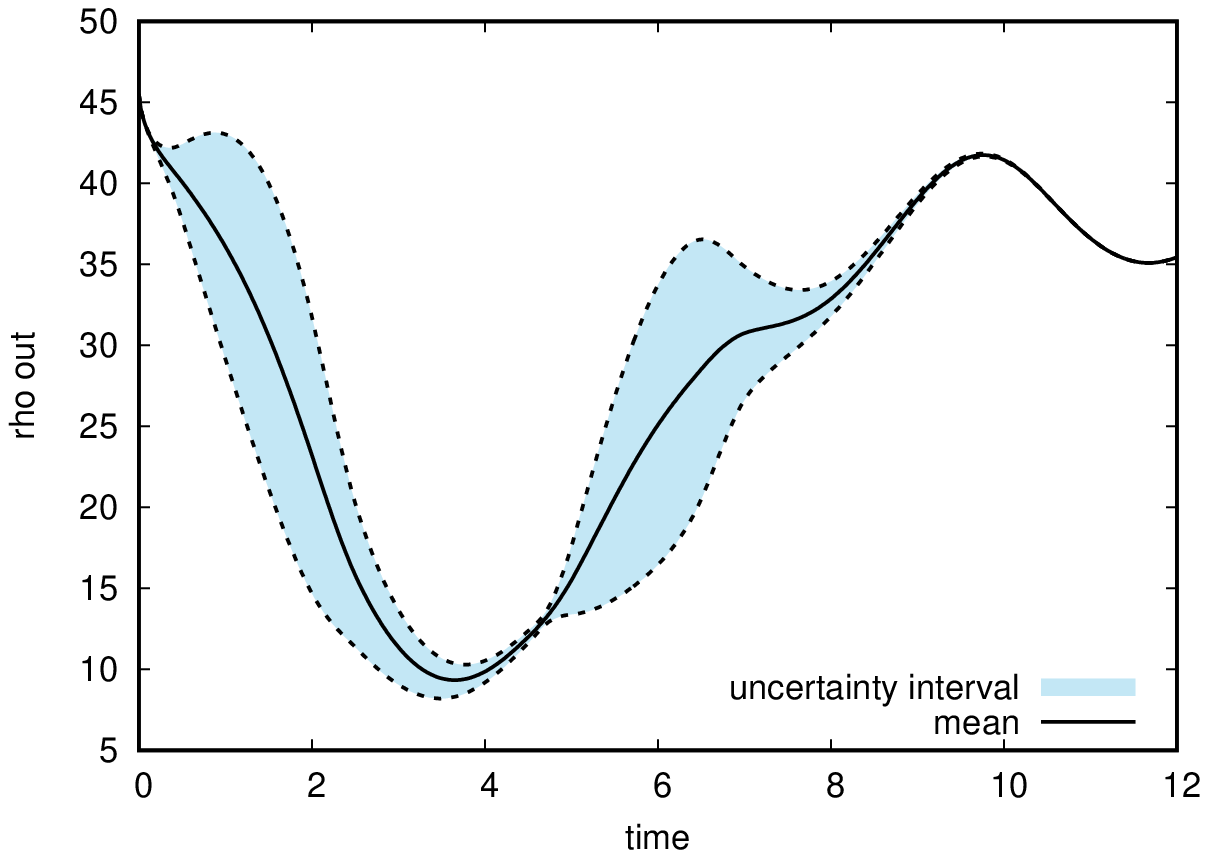} 
\caption{Uncertainty quantification results for one pipe simulation with limited inter-temporal uncertainty.  Left: inlet flow $q(t,0)$; Right: outlet density $\rho(t,L)$. The shaded region is as previously defined.}
\label{fig:comp_onepipe_3_output}
\end{figure}
In this setting, we refer to the inter-temporal uncertainty as \emph{limited}, because the perturbation may occur only on part of the simulation time interval.   
The boundary conditions are shown in Figure \ref{fig:comp_onepipe_3_bc}, where the density at the inlet is deterministic and the interval of uncertainty for the outlet flow is indicated.  The resulting solutions at the boundaries are then shown in Figure \ref{fig:comp_onepipe_3_output}. Uncertainty in inlet flow and outlet density, shown as one standard deviation above and below the mean, is seen in the at times when perturbations may occur.



\begin{figure}[h!]
\includegraphics[width=0.5\textwidth]{Fig/rho_vs_time_in_t12h_q0}   \includegraphics[width=0.5\textwidth]{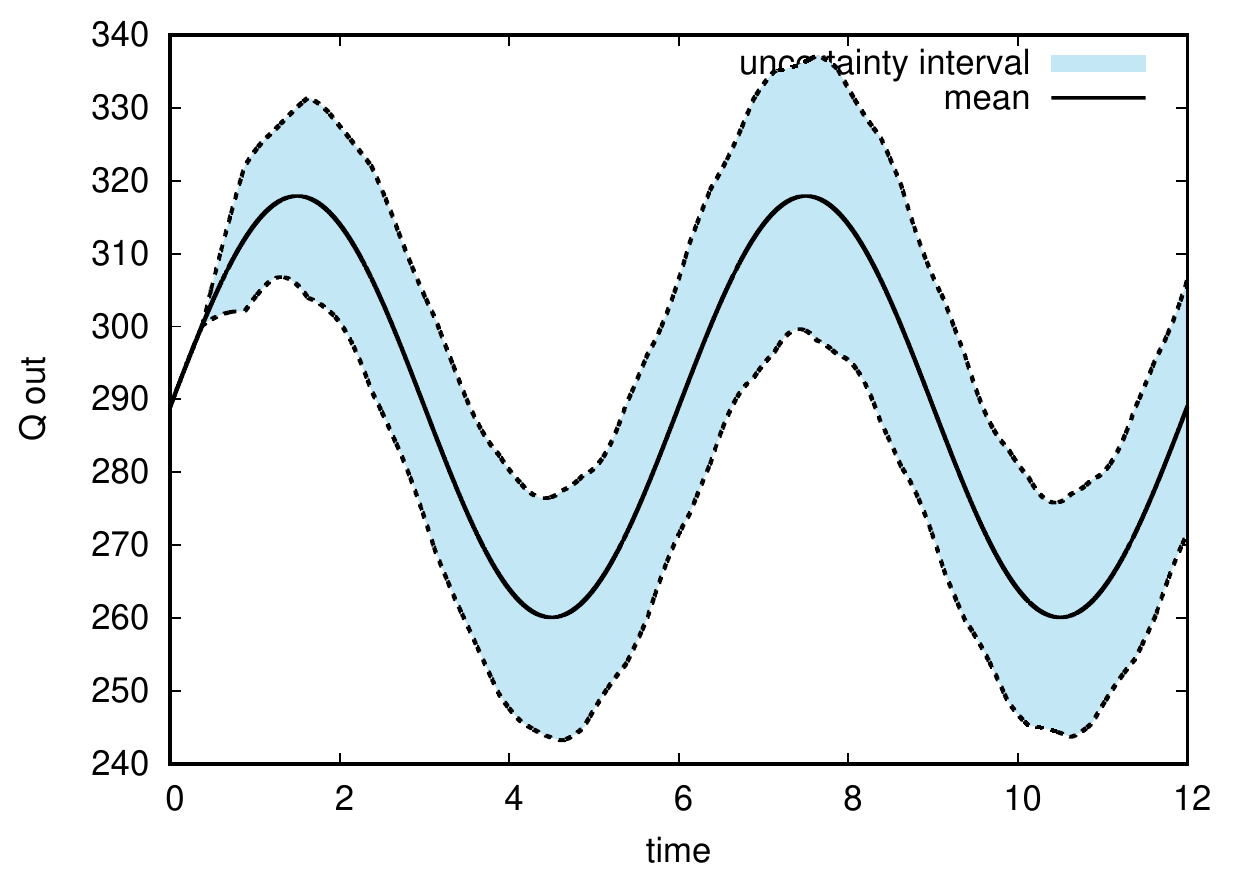} 
\caption{Boundary conditions used as inputs for one pipe simulation with limited inter-temporal uncertainty. Left: inlet density $\rho(t,0)$; Right: outlet flow $q(t,L)$.  The shaded region is as previously defined.}
\label{fig:comp_onepipe_4_bc}
\end{figure}

\begin{figure}[h!]
\includegraphics[width=0.5\textwidth]{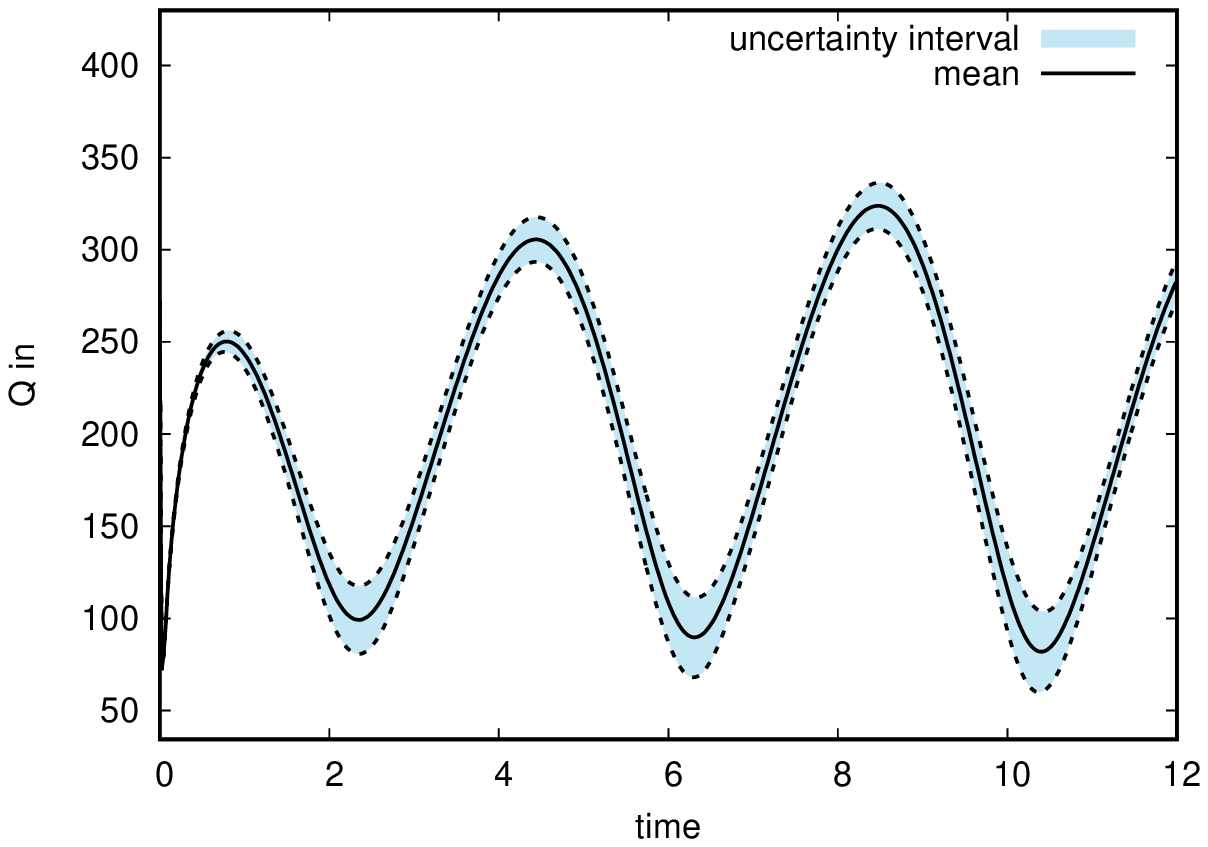}   \includegraphics[width=0.5\textwidth]{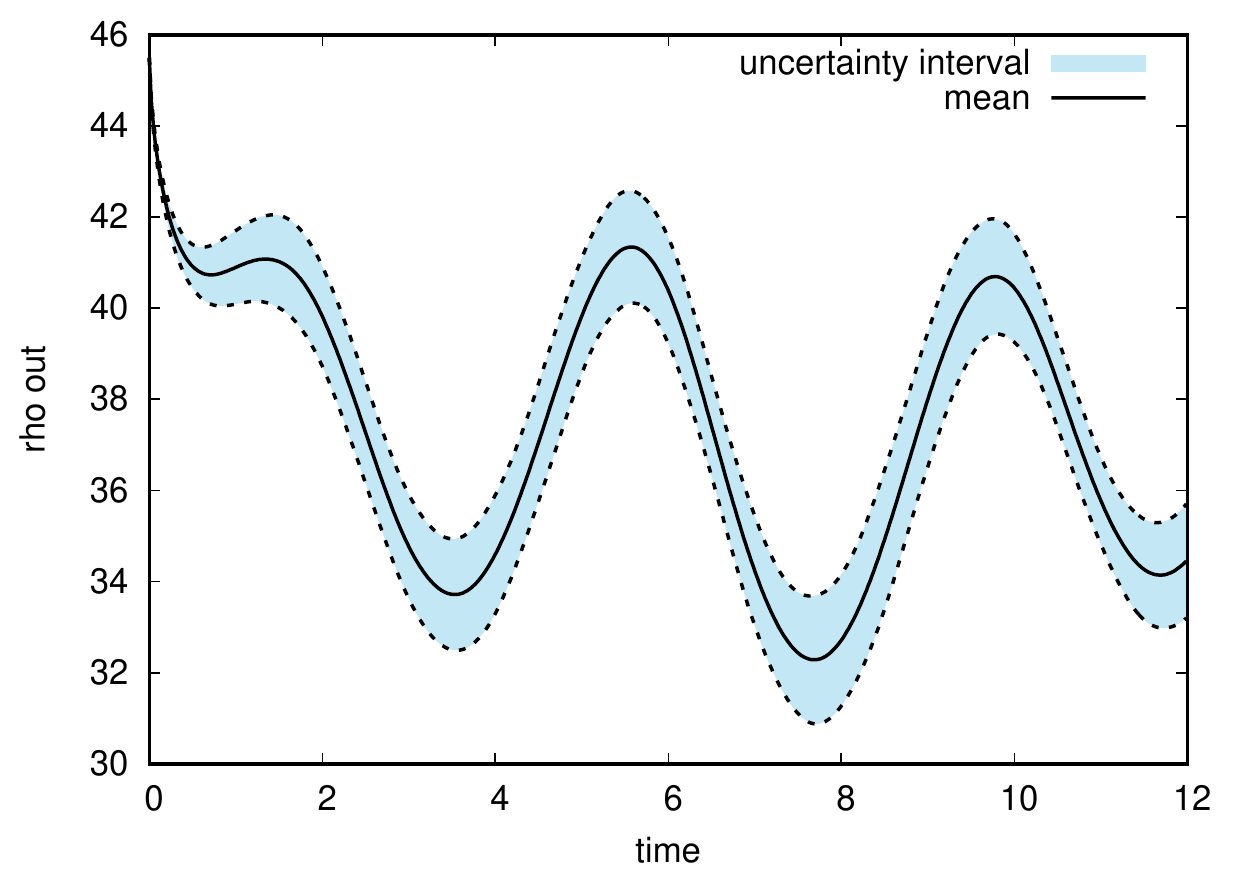} 
\caption{Uncertainty quantification results for one pipe simulation with limited inter-temporal uncertainty. Left: inlet flow $q(t,0)$; Right: outlet density $\rho(t,L)$. The shaded region is as previously defined.}
\label{fig:comp_onepipe_4_output}
\end{figure}

Alternatively, we may use $Y(\omega) \sim U[-1,11]$, so that the increased withdrawal rate can occur at any time during the simulation time interval.  We refer to this situation as \emph{universal} inter-temporal uncertainty.  For this simulation, we set $d_2(t)=1.25\cdot d_1(t)$.  The boundary conditions are shown in Figure \ref{fig:comp_onepipe_4_bc}, where the density at the inlet is deterministic and the interval of uncertainty for the outlet flow is indicated.  The resulting solutions at the boundaries are then shown in Figure \ref{fig:comp_onepipe_4_output}. Uncertainty in inlet flow and outlet density, shown as one standard deviation above and below the mean, is seen throughout the simulation interval.


          


\subsection{Uncertainty Quantification Simulations for a Test Network} \label{sec:comp_network}

\begin{figure}[h!]
\begin{center}
	\includegraphics[width=.5\linewidth]{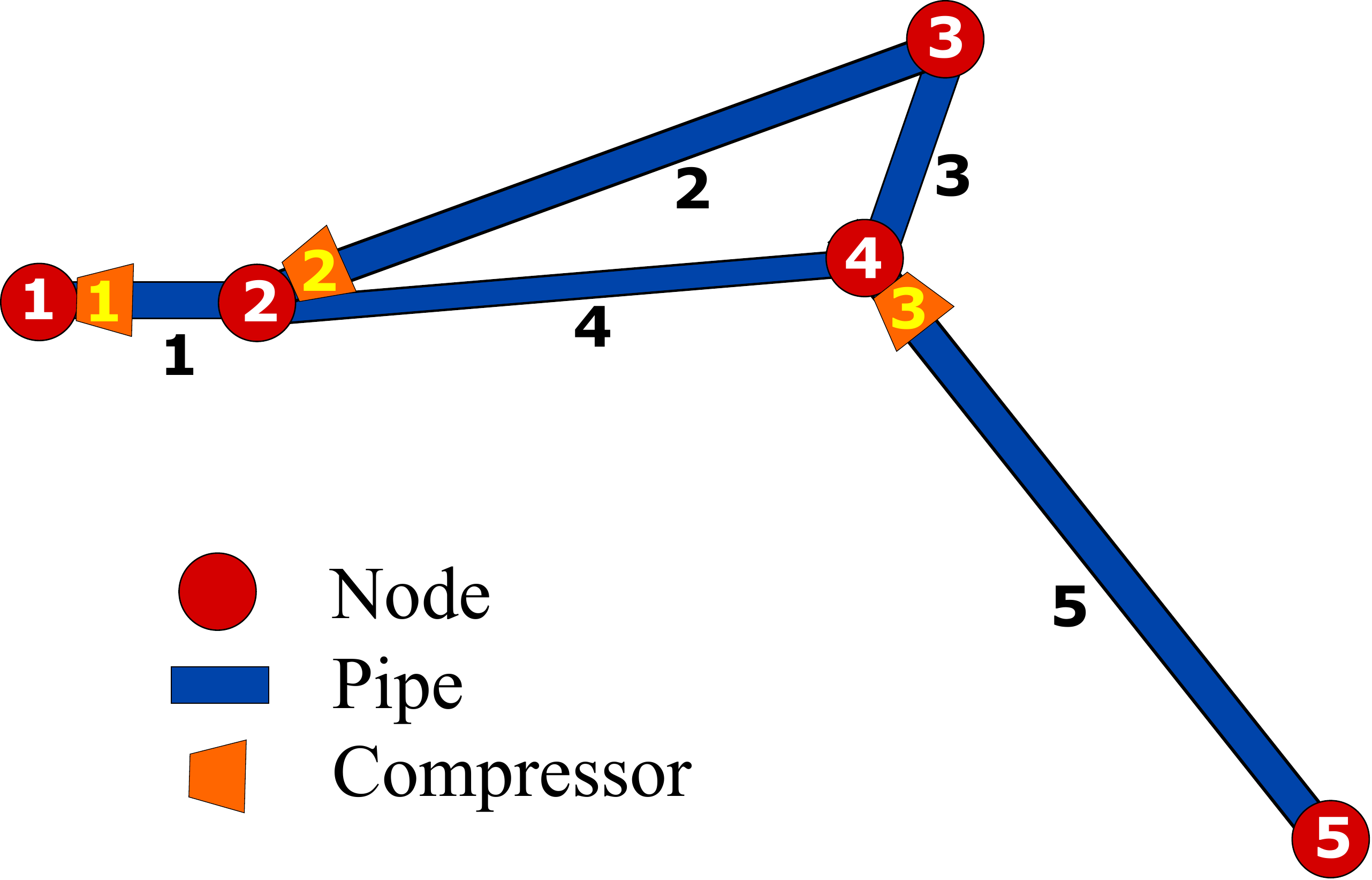}
\end{center}
\caption{5 node test network with five connecting pipes and three compressors.}
\label{fig:network}
\end{figure}

Here we simulate uncertainty propagation in a simple actuated network illustrated in Fig.~\ref{fig:network}.  The extension of the UQ analysis from the single pipe example to a network is nontrivial because of the complex Kirchoff-Neumann type boundary conditions at network nodes. The network has five nodes, five connecting pipes, and three compressors located at nodes 1, 2 and 4.  This network model has been examined in two previous pipeline simulation studies \cite{gyrya2019explicit,bermudez2021modelling}.  The network structure and random boundary conditions are fully specified here so that the stochastic test case can be reproduced, including with other methods.

The structural parameters of the network are defined for pipes in Table~\ref{tab:pipes} and for compressors in Table~\ref{tab:comps}.  This static model is defined using 5 nodes, where node 1 is a slack node with defined pressure, and nodes 2 through 5 are flow nodes with defined mass withdrawal rate.

\begin{table}[!h]
\centering
\begin{tabular}{|l|l|l|l|l|l|}
\hline
     Pipe  & From node & To node & Diameter $D_j$ (m) & Length $L_j$ (m) & Friction factor $\lambda_j$ \\ \hline
     1 & 1 & 2 & 0.9143995 & 20000 & 0.01 \\
     2 & 2 & 3 & 0.9143995 & 70000 & 0.01 \\
     3 & 3 & 4 & 0.9143995 & 10000 & 0.01 \\
     4 & 2 & 4 & 0.6349997 & 60000 & 0.015 \\
     5 & 4 & 5 & 0.9143995 & 80000 & 0.01 \\
     \hline
\end{tabular}
\caption{Pipe parameters for 5 node test network.} \label{tab:pipes}
\end{table}

\begin{table}[!h]
\centering
\begin{tabular}{|l|l|l|}
\hline
     Compressor & Location node & To pipe  \\ \hline
     1 & 1 & 1 \\
     2 & 2 & 2 \\
     3 & 4 & 5 \\
     \hline
\end{tabular}
\caption{Compressor locations for 5 node test network.} \label{tab:comps}
\end{table}

We define a stochastic initial boundary value problem as follows, using a wave speed value of $a=377.9683$ m/s as in the single pipe example described in Section \ref{sec:onepipemodel}. The initial state is deterministic and steady, and is specified using the parameters in Tables \ref{tab:initnode} and \ref{tab:initpipe}.
\begin{table}[!h]
\centering
\begin{tabular}{|l|l|l|l|}
\hline
     Item type & Item ID & Value type & Value  \\ \hline
     Node & 1 & pressure (Pa) & $\sigma_1^0=3447378.645$ \\
     Node & 2 & Flow withdrawal (kg/s) & $d_2^0=0$ \\
     Node & 3 & Flow withdrawal (kg/s) & $d_3^0=150$ \\
     Node & 4 & Flow withdrawal (kg/s) & $d_4^0=0$ \\
     Node & 5 & Flow withdrawal (kg/s) & $d_5^0=150$ \\
     Comp & 1 & Boost ratio & $c_1^0=1.5290113$ \\
     Comp & 2 & Boost ratio & $c_2^0=1.1128863$ \\
     Comp & 3 & Boost ratio & $c_3^0=1.2242249$ \\
     \hline
\end{tabular}
\caption{Initial nodal data for 5 node test network simulation.} \label{tab:initnode}
\end{table}

\begin{table}[!h]
\centering
\begin{tabular}{|l|l|l|l|}
\hline
     Pipe $k$ & Pressure in $p_k^0$ (Pa) & Pressure out $p_k^L$ (Pa) & Flow $\phi_k$ (kg/s)  \\ \hline
     1 & 5271081.1 & 4611205.3 & 300.00  \\
     2 & 5131747.2 & 3540078.3 & 233.33 \\
     3 & 3540078.3 & 3504395.3 & 83.33 \\
     4 & 4611205.3 & 3504395.3 & 66.66 \\
     5 & 4290168.0 & 3447378.6 & 150.00 \\
     \hline
\end{tabular}
\caption{Initial pipe data for 5 node test network simulation.} \label{tab:initpipe}
\end{table}

The initial conditions for pressure (in Pa) and flow (in kg/s) for each pipe $k\in\mathcal E$ are then given by
\begin{align}
\label{eq:icpipes}
p(0,x) & = \sqrt{(p_k^0)^2 - a^2\dfrac{16\lambda_k}{D_k^5\pi^2} \phi_k |\phi_k|x}, \\
\label{eq:ic-q}
\phi(0,x) & = \phi_k,
\end{align}
where $p_k^0$ and $\phi_k$ are as given in Table \ref{tab:initpipe}, and where $\lambda_k$ and $D_k$ are as given in Table \ref{tab:pipes}.  It should be straightforward to confirm that $p(0,L_k)=p_k^L$ for each pipe $k\in\mathcal E$ 

\begin{figure}[h!]
\begin{center}
	\includegraphics[width=0.99\textwidth]{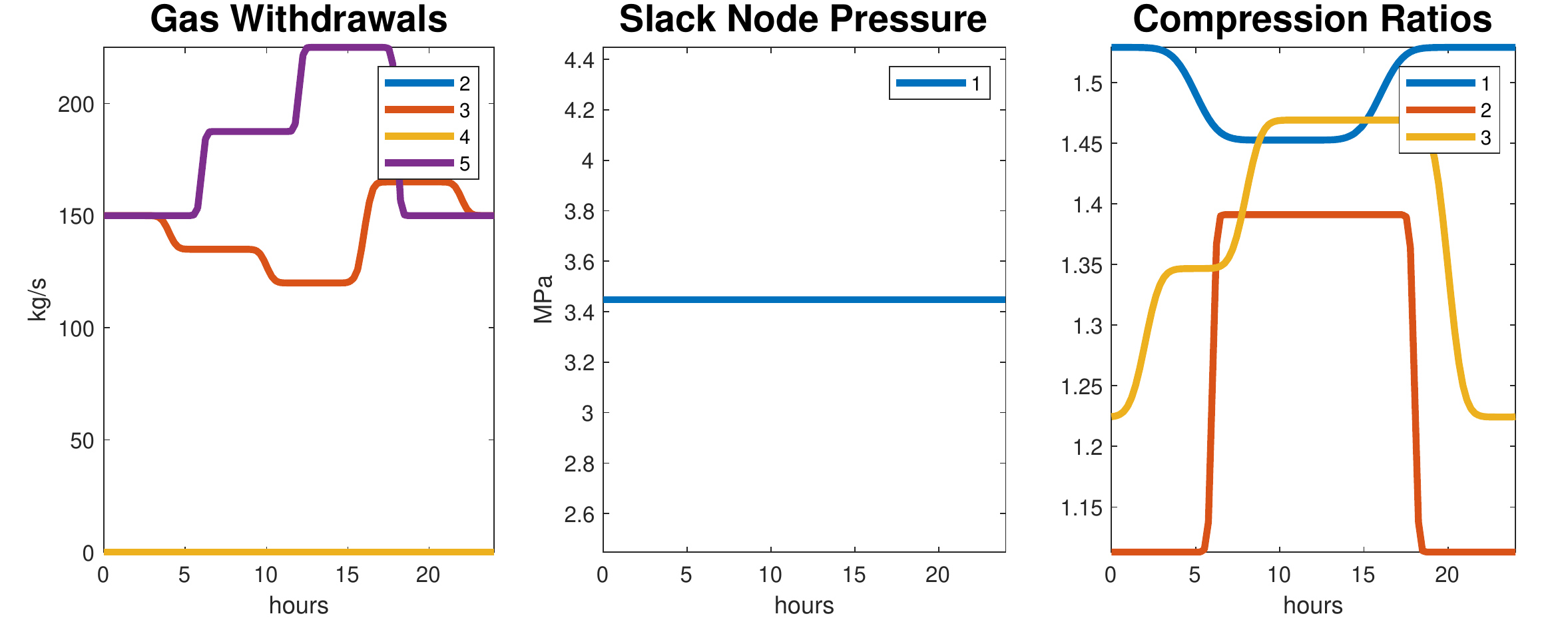}  
\end{center}
\caption{Deterministic components of boundary conditions for the 5-node test network initial boundary value problem. Left: nodal gas withdrawals $d_i(t)$ for $i=2,3,4,5$; Center: slack node pressure $\sigma_1(t)$; Right: compression ratios $c_j(t)$ for $j=1,2,3$.}
\label{fig:model5bc}
\end{figure}

The boundary conditions are defined below.  Let $d_i(t)$ denote the mass flow withdrawal (kg/s) from the network at node $i$, let $\sigma_j(t)$ denote the specified pressure at node $j$, and let $c_{k}(t)$ denote the compression ratio of compressor $k$, where time $t$ is in seconds.  The functions are defined as sums of sigmoid functions of the form $h(t)=\frac{1}{2} (\text{erf}(2t)+1)$, where $\text{erf}(t)=\frac{2}{\sqrt{\pi}}\int_0^t e^{-x^2}dx$ is the error function.
\begin{subequations}
\begin{align}
    d_3(t) & = d_3^0-d_3^0\cdot 0.1\left[h\left(\frac{t-3600\cdot 4}{3600}\right)+h\left(\frac{t-3600\cdot 10}{3600}\right)\right. \\
    & \qquad \qquad \qquad \left.-3h\left(\frac{t-3600\cdot 16}{3600}\right)+h\left(\frac{t-3600\cdot 22}{3600}\right)\right] \\
    d_5(t) & = d_5^0+d_5^0\cdot 0.25\left[h\left(\frac{t-3600\cdot 6}{1800}\right)+h\left(\frac{t-3600\cdot 12}{1800}\right)-2h\left(\frac{t-3600\cdot18}{1800}\right)\right] 
\end{align}
\end{subequations}

\begin{align}
    \sigma_1(t) & = \sigma_1^0
\end{align}

\begin{subequations}
\begin{align}
    c_1(t) & = c_1^0-c_1^0\cdot 0.05\left[h\left(\frac{t-3600\cdot 5}{1800}\right)-h\left(\frac{t-3600\cdot 16}{1800}\right)\right] \\
    c_2(t) & = c_2^0+c_2^0\cdot 0.25\left[h\left(\frac{t-3600\cdot 6}{1800}\right)-h\left(\frac{t-3600\cdot 18}{1800}\right)\right] \\
    c_3(t) & = c_3^0+c_3^0\cdot 0.1\left[h\left(\frac{t-3600\cdot 2}{7200}\right)+h\left(\frac{t-3600\cdot 8}{7200}\right)-2 h\left(\frac{t-3600\cdot 20}{7200}\right)\right]
\end{align}
\end{subequations}

\begin{figure}[h!]
\begin{center}
	\includegraphics[width=0.49\textwidth]{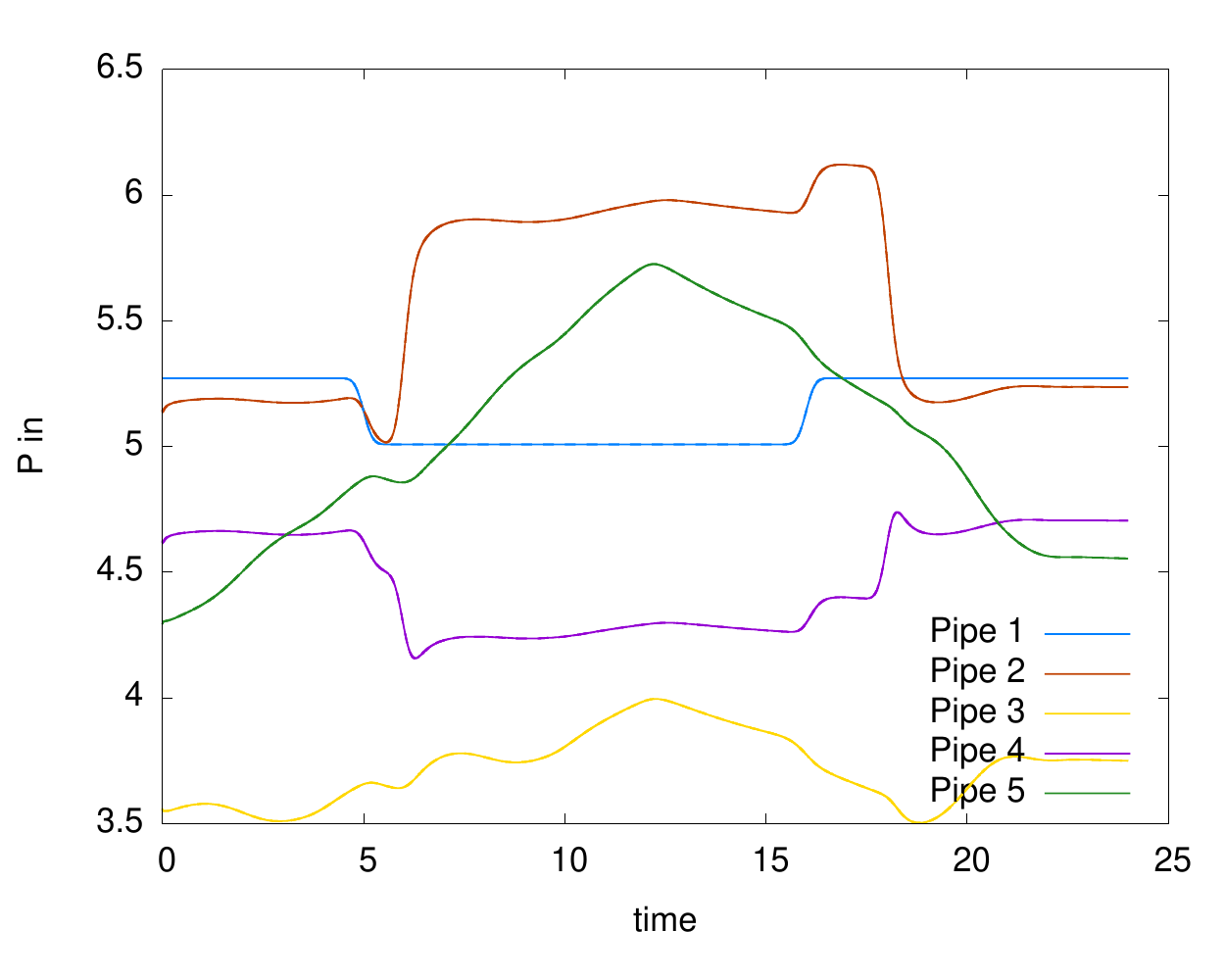}   \includegraphics[width=0.49\textwidth]{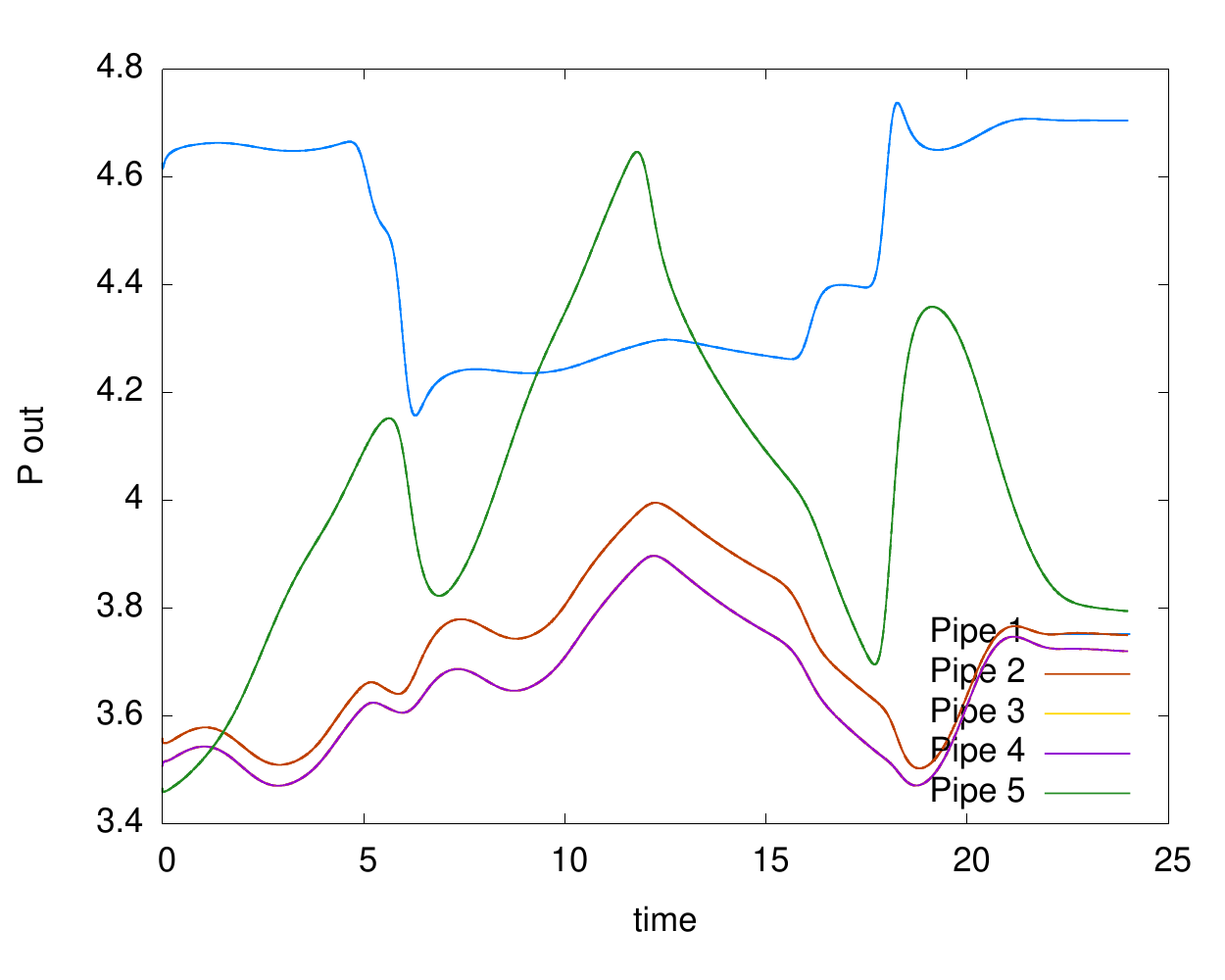} \\
	\includegraphics[width=0.49\textwidth]{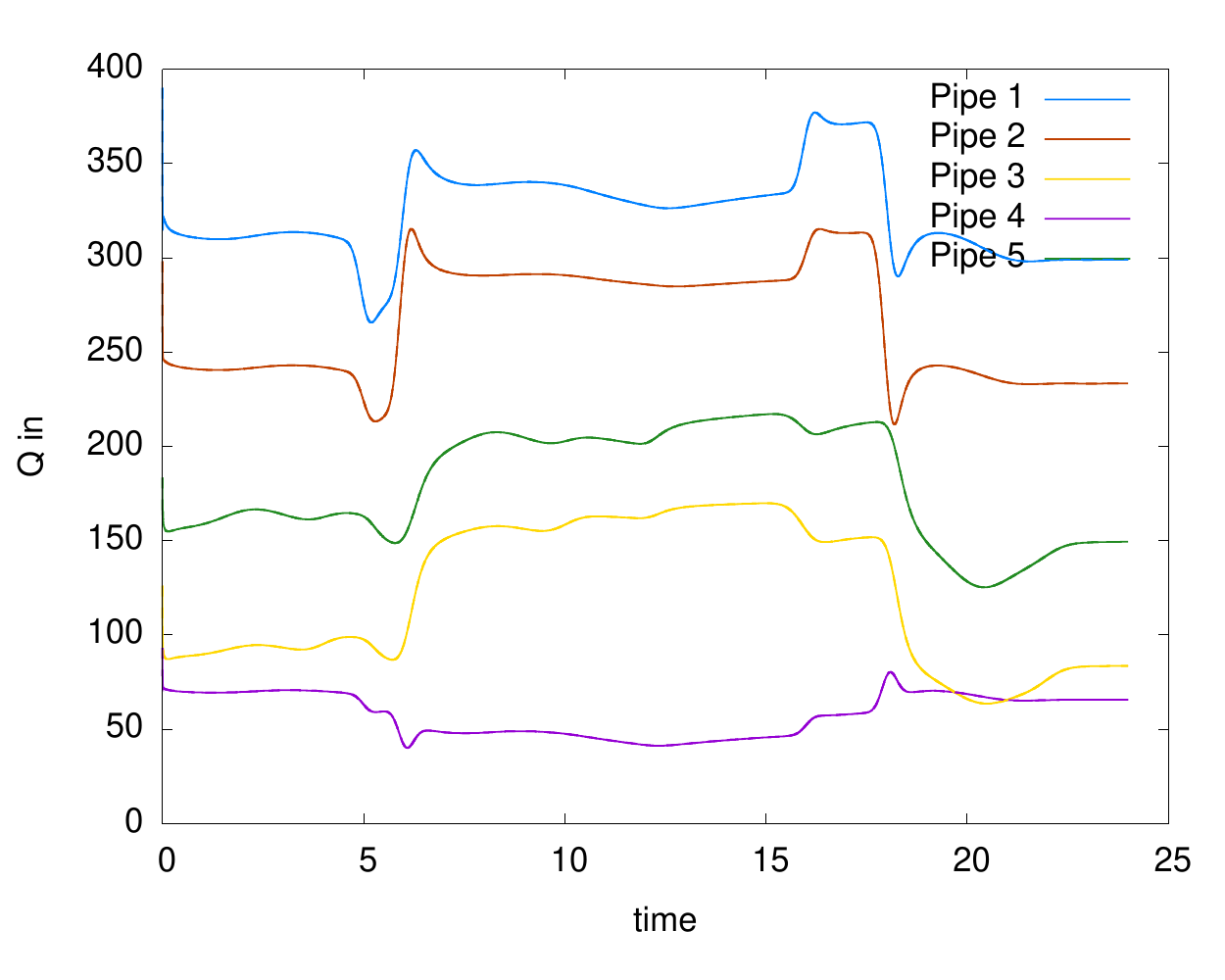}   \includegraphics[width=0.49\textwidth]{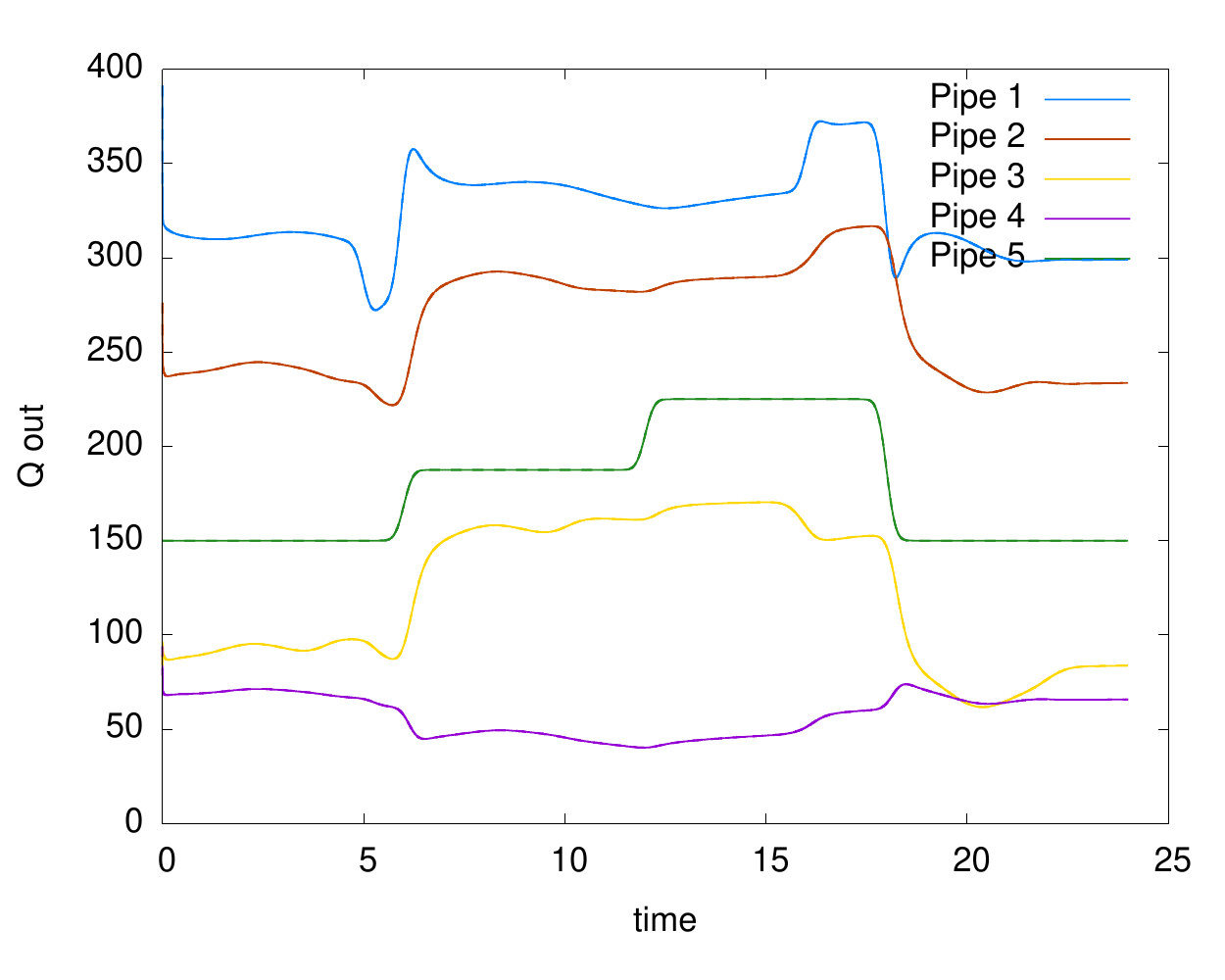} 
\end{center}
\caption{Deterministic solution for the 5-node test network initial boundary value problem using boundary conditions defined in Figure \ref{fig:model5bc}. Top: inlet (left) and outlet (right) pressure (MPa) for each pipe. Bottom: inlet (left) and outlet (right) flow (kg/s) for each pipe.}
\label{fig:model5detsol}
\end{figure}

These deterministic components of the boundary conditions are shown in Figure \ref{fig:model5bc}, and the resulting deterministic solutions for pressure and flow at the pipe endpoints are shown in Figure \ref{fig:model5detsol}.  The deterministic problem is then extended to the stochastic setting by adding intertemporal uncertainty to a single boundary parameter.

\begin{figure}[h!]
\begin{center}
	\includegraphics[width=0.49\textwidth]{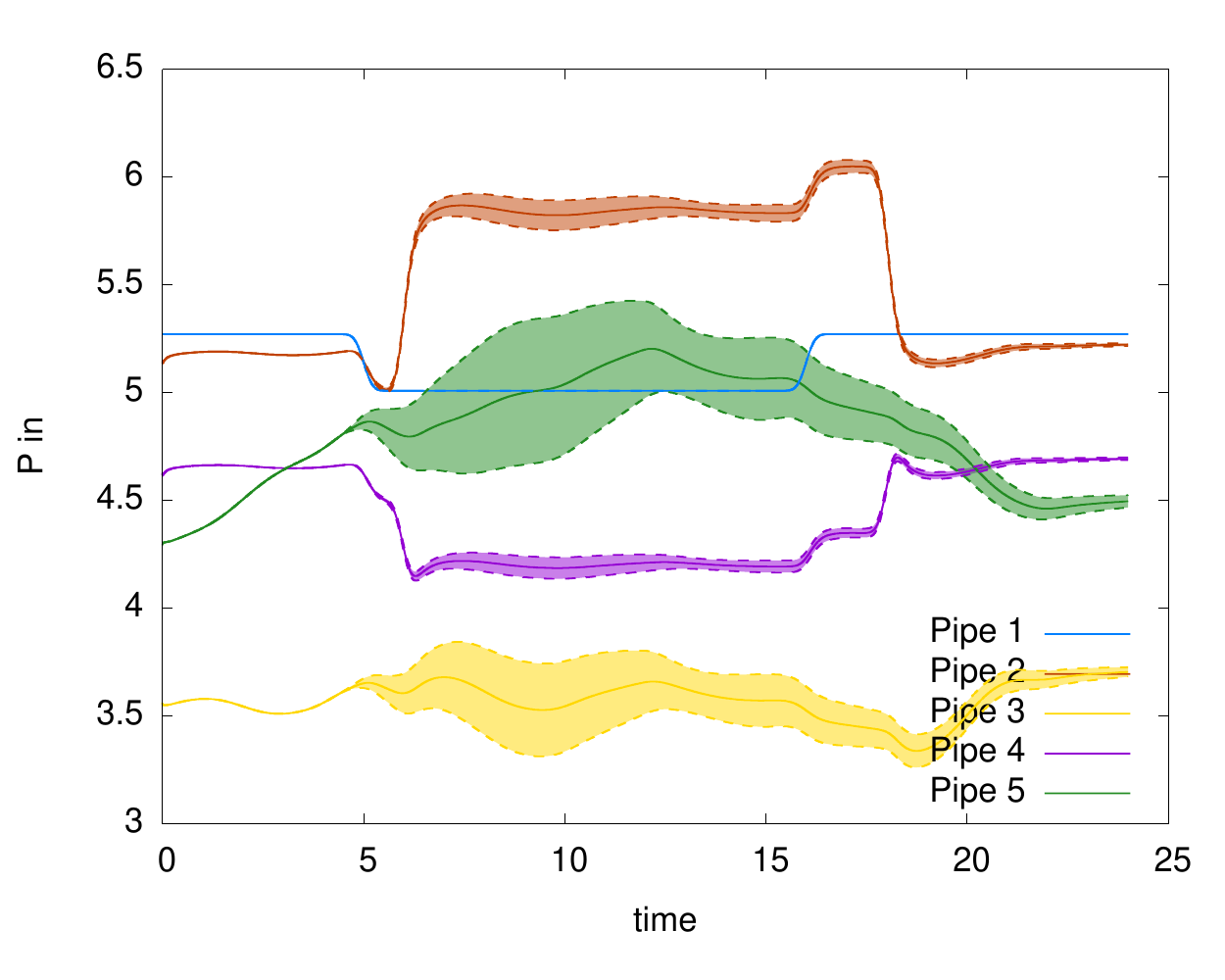}   \includegraphics[width=0.49\textwidth]{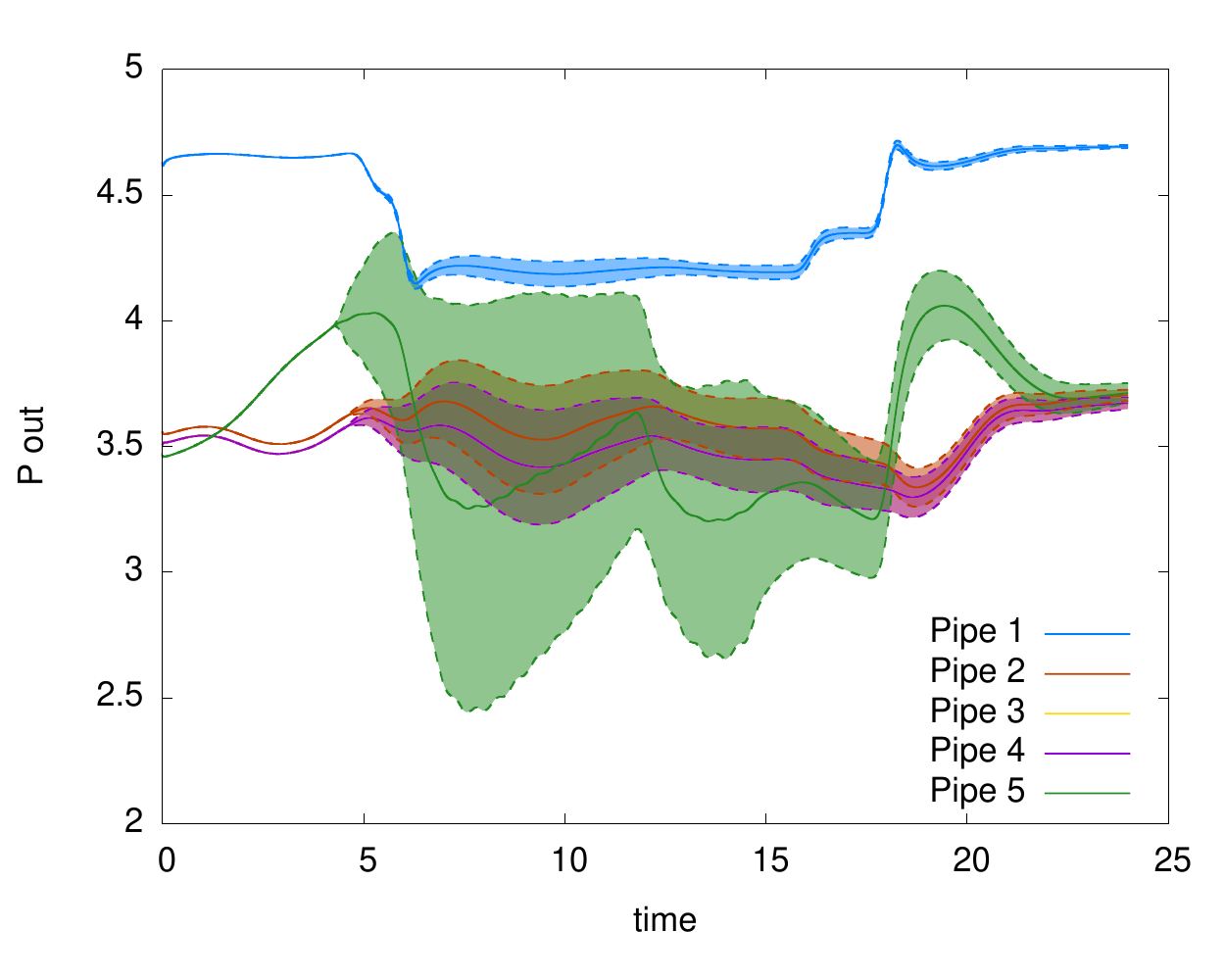} \\
	\includegraphics[width=0.49\textwidth]{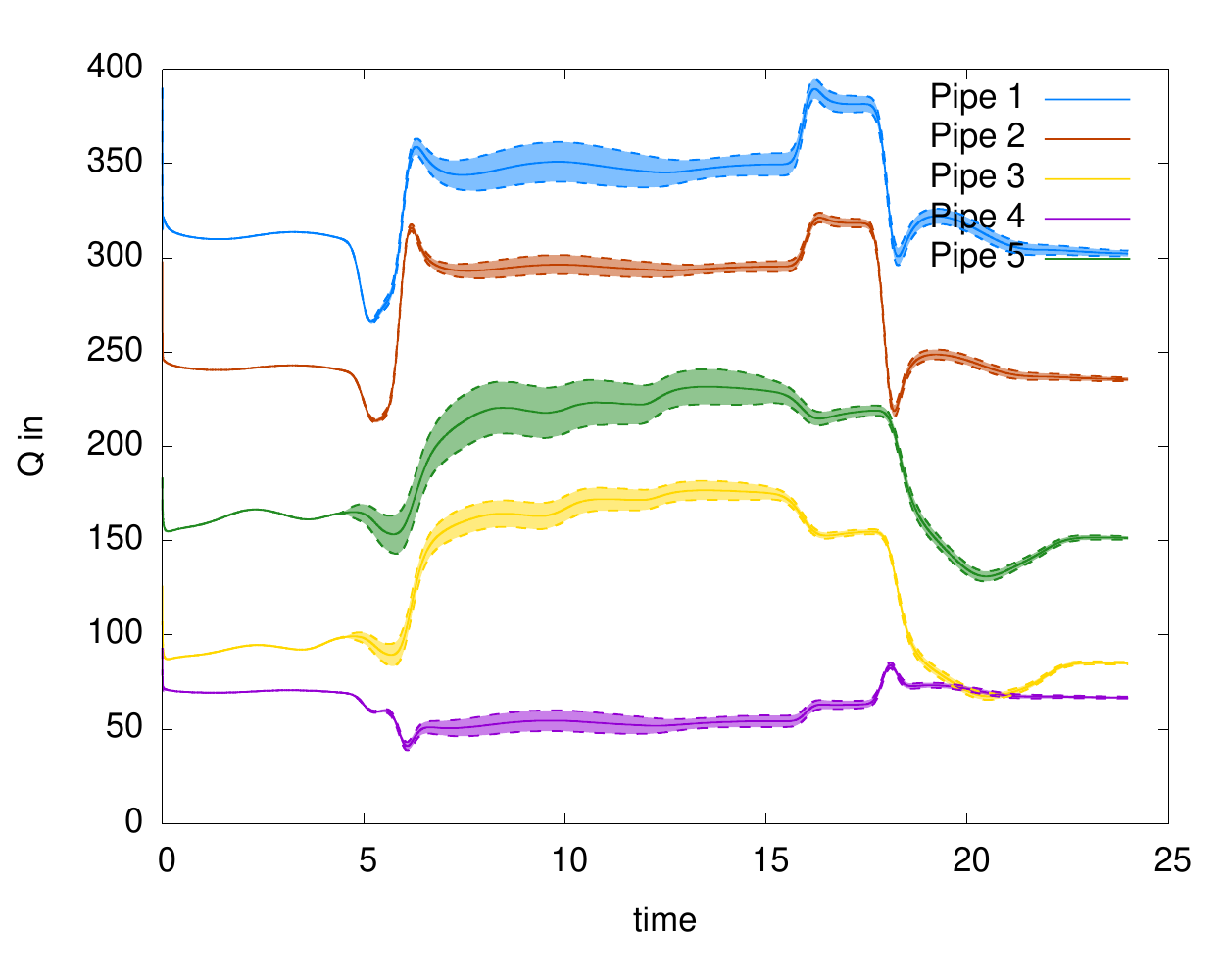}   \includegraphics[width=0.49\textwidth]{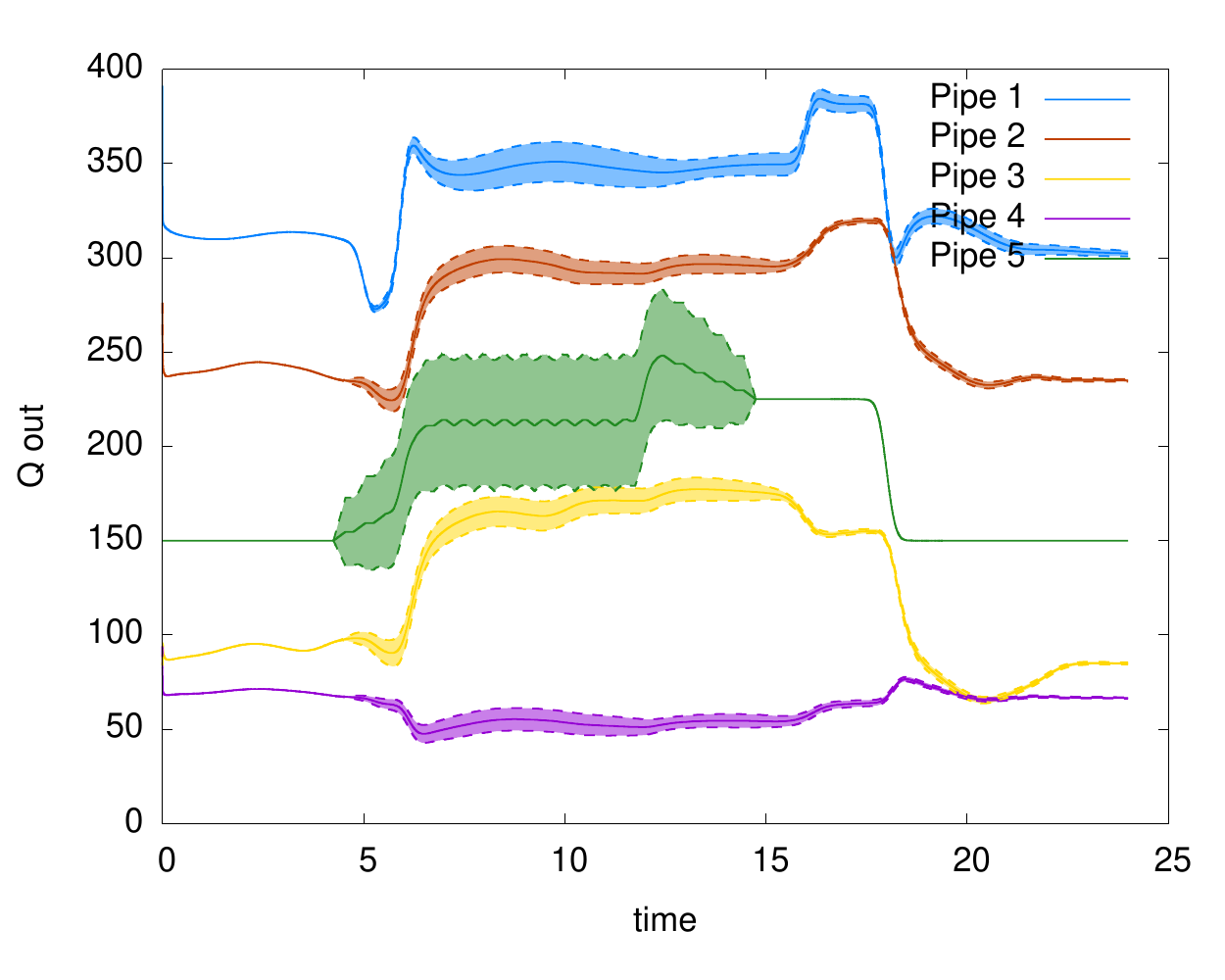} 
\end{center}
\caption{Uncertainty quantification results for the 5-node test network. Top: inlet (left) and outlet (right) pressure (MPa) for each pipe. Bottom: inlet (left) and outlet (right) flow (kg/s) for each pipe.}
\label{fig:network-results-uq}
\end{figure}

To create a stochastic IBVP, we add intertemporal uncertainty to the withdrawal rate at node $5$ in the same way as described for a single pipe in the example in Section \ref{sec:example3}. For this simulation, we have chosen $N_y = 16$ cells in the stochastic direction to ensure an accurate representation of the input uncertainty at the outflow boundary of the 5-th pipe. The results of uncertainty quantification using the SFV method for $T_p(\omega) \sim U[4,12]$ and $d_2  = d_1 + 0.5d_5^0$ are shown in Fig.~\ref{fig:network-results-uq}. As in the single pipe examples, a shaded region that extends one standard deviation above and below the mean is superimposed on each time-varying solution to indicate uncertainty. We clearly observe that the intertemporal gas consumption uncertainty at node 5 propagates quickly through the network.

\begin{figure}[h!]
    \centering
    \includegraphics[width=0.49\textwidth]{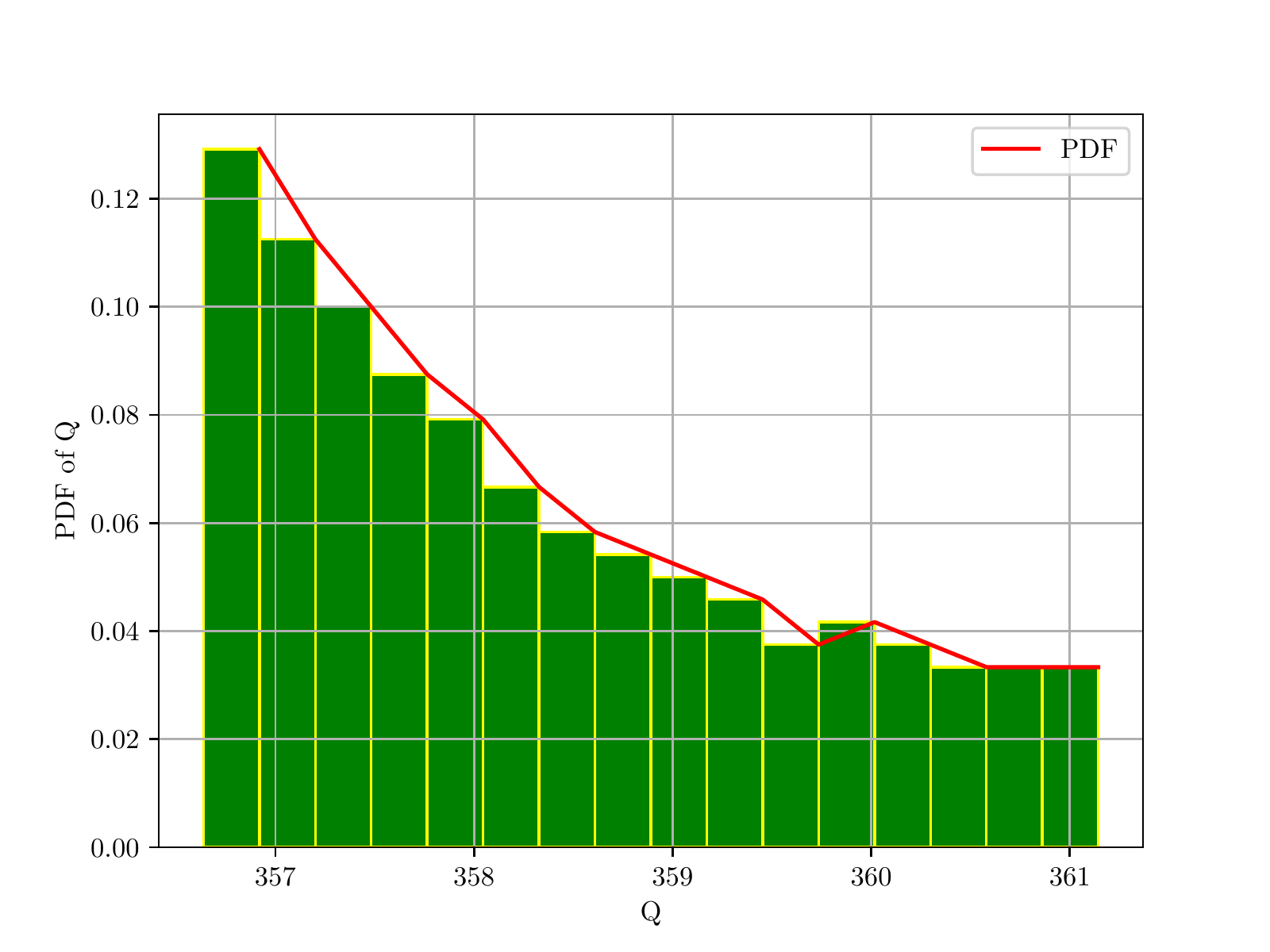} \includegraphics[width=0.49\textwidth]{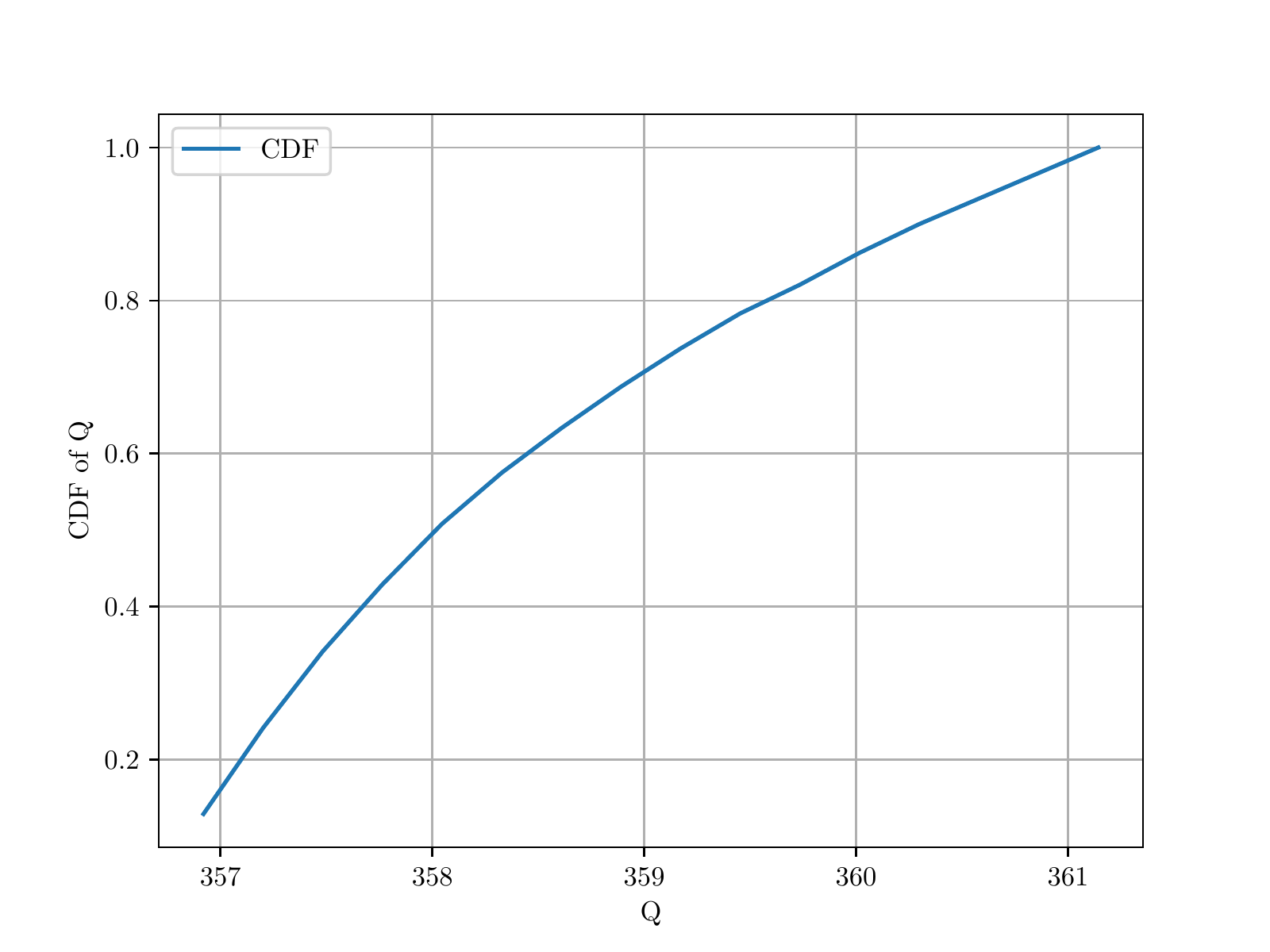}
    \caption{Probability density function (left) and cumulative distribution function (right) for $Q$ in the middle of the pipe $2$.}
    \label{fig:distr}
\end{figure}

In addition, a rapid computation can be applied to the solution produced by the SFV method to extract the probability density functions of the quantities of interest at any stage of the simulation. For example, Fig.~\ref{fig:distr} shows the probability density function (left) and cumulative distribution function (right) for flow $q_2(12\cdot 3600,L_2/2)$ in the middle of the pipe 2 at hour 12 of the simulation.  The distribution is non-trivially non-symmetric, and reflects the nonlinear nature of uncertainty propagation through the network of hyperbolic conservation laws.

\section{Discussion}

The SFV method that we present here is a computationally efficient technique for fully characterizing a stochastic process that arises from an IBVP for hyperbolic conservation laws on domains with boundaries connected in a graph.  The approach can account for uncertainty in initial conditions, model parameters, and complex time-varying random boundary conditions.  We have demonstrated through computational experiments (see Figure \ref{fig:onepipe_convergence}) that the accuracy of the scheme improves with increasing approximation order, so that the method is significantly less computationally costly than Monte Carlo simulation.  Moreover, simple computations in the stochastic domain can be done to obtain statistics of interest at a time-slice such as moments, or the distribution itself.  

The results of the computational method are somewhat sensitive to the choice of discretizations in space, time, and stochastic domains.  In particular, if discretization in the stochastic space is too coarse, this can lead to spurious time-dependent artifacts in the evaluation of inter-temporal uncertainty propagation.  The discretization should be sufficiently fine to appropriately resolve the solution given the parameters of the uncertainty.  Guidelines for choices of such computational parameters, as well as theoretical underpinnings related to accuracy and convergence, are open questions that would be of interest in future studies. 

Our study on efficient uncertainty quantification for hyperbolic conservation laws is motivated by the need to model and further to optimize large-scale natural gas pipeline flows in a model-predictive manner in the presence of uncertainty.  Here, we enable the efficient characterization of both interval and inter-temporal uncertainty, which can be used to determine the capabilities of a natural gas transmission system to accommodate variable and intermittent consumption.  Precise quantification of the statistical distributions in pressures and flows can then be used to provide probabilistic guarantees on system integrity and the temporary reserves of energy that a pipeline stores in the form of ``line-pack'', or mass of gas in the pipe in the neighborhood of a consumption point.  The efficiency and flexibility of the SFV method makes it promising for use in iterative methods that could calibrate parameters such as compressor settings and maximal throughput.
	
\section{Conclusions}

We have developed a novel Stochastic Finite Volume method for stochastic hyperbolic PDEs on graphs and demonstrated the robustness and accuracy of the method for uncertainty quantification of gas flow on networks. The method is semi-intrusive, therefore only minor changes to the underlying deterministic initial boundary value problem solver are required. The solution obtained by the SFV method yields a complete solution in physical space and in time, as well as in the stochastic space.  Therefore, all  stochastic information of interest can be easily extracted, including confidence intervals, distributions, and moments. This motivates potential extensions of this method as a promising tool with which to develop methods for robust optimization of hyperbolic PDE systems under uncertainty.

\bibliographystyle{unsrt}
\bibliography{biblio,references}



\end{document}